\documentclass[12pt]{article}
\usepackage{latexsym,amsfonts,amssymb}
\usepackage[dvips]{color}
\usepackage{colordvi,multicol}
\usepackage{amsmath}
\usepackage{graphicx}
\usepackage{pdflscape}
\usepackage{rotating}

\usepackage{latexsym,amsfonts,amssymb}
\usepackage[dvips]{color}
\usepackage{colordvi,multicol}
\usepackage{amsmath}
\usepackage{graphicx}
\usepackage{pdflscape}
\usepackage{rotating}
\usepackage{breqn}

\def \cal{\mathcal}

\newtheorem{thm}{Theorem}[section]

\newtheorem{lem}[thm]{Lemma}
\newtheorem{pro}[thm]{Proposition}

\newtheorem{rem}[thm]{Remark}

\newtheorem{que}[thm]{Question}

\date{}
\begin{document}

\title{\bf An equivalent inequality for the Riemann
hypothesis}
\author{Wei Sun\\ \\ \\
  {\small Department of Mathematics and Statistics}\\
{\small Concordia University}\\
{\small Montreal, H3G 1M8,  Canada}\\
{\small  wei.sun@concordia.ca}}

\maketitle

\begin{abstract}
\noindent We present a purely analytical inequality which is equivalent to the Riemann hypothesis (RH). The proof of the equivalence is based on a representation of the modulus of the Riemann $\xi$ function. As the first step to analyze the inequality, we consider polynomial approximations. We also show that the RH is equivalent to the statement that some wave functions
constructed using the Brownian motion never evolve into perfectly distinguishable states.
\end{abstract}

\noindent  {\it MSC:} 11M26, 26D05, 26D15, 81Q10.

\noindent  {\it Keywords:} Riemann hypothesis, inequality, Riemann $\xi$ function, characteristic function, polynomial approximation, wave function, orthogonalization time.

\section{Introduction}\setcounter{equation}{0}

The Riemann hypothesis (RH) is the conjecture that the Riemann zeta function $\zeta(s)$ has its nontrivial zeros only on the critical line ${\rm Re}(s)=\frac{1}{2}$. Let $\xi(s)$ be the Riemann xi function, i.e.,
$$
\xi(s)=\frac{1}{2}s(s-1)\pi^{-\frac{s}{2}}\Gamma\left(\frac{s}{2}\right)\zeta(s),\ \ \ \ s\in\mathbb{C},
$$
where $\Gamma(s)$ is the Gamma function. Then, the RH is equivalent to the statement that the zeros of $\xi(s)$  are all located on the critical line. We refer the reader to Edwards's book \cite{E} for the terminology used in this paper, and to Wikipedia \cite{W} and references therein for an introduction to the RH.

For $\tau\in\mathbb{R}$, we define the function ${\cal J}_{\tau}$ by
$$
{\cal J}_{\tau}(y)=\sum_{m=1}^{\infty}\sum_{n=1}^{\infty}\left(\frac{n}{m}\right)^{\tau}e^{-2\pi mny},\ \ \ \ y>0,
$$
and
the function ${\eta}_{\tau}$ by
$$
{\eta}_{\tau}(y)=\frac{(y+\sqrt{y^2-1})^{\tau}+(y+\sqrt{y^2-1})^{-\tau}}{\sqrt{y^2-1}},\ \ \ \ y>1.
$$
In this paper, we will prove the following result.
\begin{thm}\label{thm0000}
The RH is equivalent to the statement that for any $\tau\in(0,\frac{1}{2})$,
\begin{eqnarray}\label{Jan27a}
0&<&4\left[\left(t^2+\tau^2+\frac{1}{4}\right)^2-\tau^2\right]\int_1^{\infty}\int_1^{\infty}\cos(2t\ln x){\cal J}_{\tau}(x^2y){\eta}_{\tau}(y)dxdy\nonumber\\
&&+\left(t^2+2\tau^2+\frac{1}{4}\right)\int_1^{\infty}\left[2y{\cal J}'_{\tau}(y)+{\cal J}_{\tau}(y)\right]{\eta}_{\tau}(y)dy\nonumber\\
&&-\int_1^{\infty}y\left[2y^2{\cal J}'''_{\tau}(y)+9y{\cal J}''_{\tau}(y)+6{\cal J}'_{\tau}(y)\right]{\eta}_{\tau}(y)dy,\ \ \ \ \ \ \ \ \forall t\in\mathbb{R}.
\end{eqnarray}
\end{thm}

The rest of the paper is organized as follows. In Section 2, we first make some preparation. Then, we derive in Section 3 a  modulus representation of the function $\xi$ (see Theorem \ref{thm000}). Theorem \ref{thm0000} is a direct consequence of the modulus representation. In Section 4, we discuss polynomial approximations of a variant of inequality (\ref{Jan27a}). In Section 5, we present another variant of inequality (\ref{Jan27a}) (see Theorem \ref{thm1}) and connect the RH with the orthogonalization time of quantum states. Through considering a  quantum dynamics defined on the classical Wiener space and constructing suitable wave functions
by virtue of the Brownian motion, we will show that the RH is equivalent to the statement that the wave functions never evolve into orthogonal states (see Theorem \ref{thm2}). Some numerical calculations for the key constants of the variant inequalities are given in the Appendix.

\section{Preliminary}\setcounter{equation}{0}

For $\sigma\in\mathbb{R}$, we define
$$
\Xi_{\sigma}(t)=\frac{\xi(\sigma-it)}{\xi(\sigma)},\ \ \ \ t\in\mathbb{R}.
$$
\begin{thm}\label{thm00} (\cite[Nakamura]{N}) The function $\Xi_{\sigma}(t)$ is a characteristic function for any $\sigma\in\mathbb{R}$. Moreover,
the probability density function $P_{\sigma}(y)$ is given as follows:
$$
P_{\sigma}(y)= \left\{ \begin{array}{ll}
        \frac{1}{\xi(\sigma)}\sum\limits_{n=1}^{\infty}f(ne^{-y})e^{-\sigma y}, \ \ & y\le 0,\\
        \frac{1}{\xi(\sigma)}\sum\limits_{n=1}^{\infty}f(ne^{y})e^{(1-\sigma)y}, \ \ & y> 0,\end{array} \right.
$$
where $f(x)=2\pi(2\pi x^4-3x^2)e^{-\pi x^2}$.

\end{thm}

\begin{rem}
In \cite[Theorem 1.1]{N}, $\frac{2}{\xi(\sigma)}$ instead of $\frac{1}{\xi(\sigma)}$ is used. But the function $\xi(s)$ used in \cite{N} is equal to $2\xi(s)$ of our paper (cf.  \cite[(1.1)]{N}). So $\frac{2}{\xi(\sigma)}$ should be replaced by $\frac{1}{\xi(\sigma)}$ here according to our notation.
\end{rem}

For $y>0$, define
\begin{eqnarray}\label{Jan27b1}
R(y)=2\sum_{n=1}^{\infty}e^{-\pi n^2y^2},
\end{eqnarray}
and
$$
G(y)=1+R(y),\ \ \ \ H(y)=y[yG(y)]''.
$$
We have
\begin{eqnarray*}
H(y)=y^2G''(y)+2yG'(y)=4\pi\sum\limits_{n=1}^{\infty}[2\pi(ny)^4-3(ny)^2]e^{-\pi n^2y^2},\ \ \ \ y>0,
\end{eqnarray*}
and
\begin{eqnarray}\label{12/11}
P_{\sigma}(y)= \left\{ \begin{array}{ll}
        \frac{1}{2\xi(\sigma)}H(e^{-y})e^{-\sigma y}, \ \ & y\le 0,\\
        \frac{1}{2\xi(\sigma)}H(e^{y})e^{(1-\sigma)y}, \ \ & y> 0.\end{array} \right.
\end{eqnarray}

Let $Y_{\sigma}$ and  $Y'_{\sigma}$ be two independent continuous random variables with the same probability density function $P_{\sigma}(y)$. Denote by $\bar{P}_{\sigma}(y)$ the probability density function of $Y_{\sigma}-Y'_{\sigma}$. By (\ref{12/11}), for $y\ge0$, we have that
\begin{eqnarray*}
\bar{P}_{\sigma}(y)
&=&\int_{-\infty}^{\infty}P_{\sigma}(y-z)P_{\sigma}(-z)dz\\
&=&\int_{-\infty}^{\infty}P_{\sigma}(y+z)P_{\sigma}(z)dz\\
&=&\int_{0}^{\infty}P_{\sigma}(y+z)P_{\sigma}(z)dz+\int_{-y}^{0}P_{\sigma}(y+z)P_{\sigma}(z)dz+\int_{-\infty}^{-y}P_{\sigma}(y+z)P_{\sigma}(z)dz\\
&=&\frac{1}{4\xi^2(\sigma)}\Bigg[\int_{0}^{\infty}H(e^{y+z})e^{(1-\sigma)(y+z)}H(e^{z})e^{(1-\sigma) z}dz\\
&&\ \ \ \ +\int_{-y}^{0}H(e^{y+z})e^{(1-\sigma)(y+z)}H(e^{-z})e^{-\sigma z}dz+\int_{-\infty}^{-y}H(e^{-(y+z)})e^{-\sigma(y+z)}H(e^{-z})e^{-\sigma z}dz\Bigg]\\
&=&\frac{1}{4\xi^2(\sigma)}\Bigg[e^{(1-\sigma)y}\int_{0}^{\infty}H(e^{y+z})H(e^z)e^{2(1-\sigma)z}dz+e^{(1-\sigma)y}\int_{0}^{y}H(e^{y-z})H(e^z)e^{(2\sigma-1)z}dz\nonumber\\
&&\ \ \ \ +e^{\sigma y}\int_{0}^{\infty}H(e^{y+z})H(e^z)e^{2\sigma z}dz\Bigg].
\end{eqnarray*}

For $y\ge0$, define
\begin{eqnarray}\label{Oct22a}
{\cal U}_{\sigma}(y)&=&e^{(1-\sigma)y}\int_{0}^{\infty}H(e^{y+z})H(e^z)e^{2(1-\sigma)z}dz+e^{(1-\sigma)y}\int_{0}^{y}H(e^{y-z})H(e^z)e^{(2\sigma-1)z}dz\nonumber\\
&&+e^{\sigma y}\int_{0}^{\infty}H(e^{y+z})H(e^z)e^{2\sigma z}dz.
\end{eqnarray}
Then,
\begin{eqnarray*}
\bar{P}_{\sigma}(y)
=\frac{{\cal U}_{\sigma}(|y|)}{4\xi^2(\sigma)},\ \ \ \ y\in\mathbb{R}.
\end{eqnarray*}
Hence,
\begin{eqnarray}\label{RRR2}
|\xi(\sigma-it)|^2=\xi^2(\sigma)E\left[e^{it(Y_{\sigma}-Y'_{\sigma})}\right]=\frac{1}{2}\int_0^{\infty}{\cal U}_{\sigma}(y)\cos(ty)dy.
\end{eqnarray}

\begin{rem} By (\ref{Oct22a}), for $\sigma\in(\frac{1}{2},1)$ and $y\ge 0$, we get
\begin{eqnarray*}
\frac{{\cal U}_{\sigma}(y)}{4\pi}&=&\sum\limits_{n=1}^{\infty}\Bigg\{e^{(1-\sigma)y}\int_{0}^{\infty}[2\pi (ne^{y+z})^4-3(ne^{y+z})^2]e^{-\pi (ne^{y+z})^2}H(e^z)e^{2(1-\sigma)z}dz\\
&&\ \ \ \ \ \ \ +e^{(1-\sigma)y}\int_{0}^{y}[2\pi (ne^{y-z})^4-3(ne^{y-z})^2]e^{-\pi (ne^{y-z})^2}H(e^z)e^{(2\sigma-1)z}dz\\
&&\ \ \ \ \ \ \ +e^{\sigma y}\int_{0}^{\infty}[2\pi (ne^{y+z})^4-3(ne^{y+z})^2]e^{-\pi (ne^{y+z})^2}H(e^z)e^{2\sigma z}dz\Bigg\}\\
&<&2\pi e^{\frac{9}{2}y}\sum\limits_{n=1}^{\infty}n^4\int_0^{\infty}e^{5z}\left[e^{-\pi n^2e^{2z}}\right]^{e^{2y}}{H}(e^z)dz\nonumber\\
&&+2\pi e^{\frac{9}{2}y}\sum\limits_{n=1}^{\infty}n^4e^{-\pi n^2e^{y}}\int_{0}^{\frac{y}{2}}{H}(e^z)dz+2\pi e^{3y}\sum\limits_{n=1}^{\infty}n^4e^{-\pi n^2}\int_{\frac{y}{2}}^{y}{H}(e^z)dz\nonumber\\
&&+2\pi e^{5y}\sum\limits_{n=1}^{\infty}n^4\int_0^{\infty}e^{6z}\left[e^{-\pi n^2e^{2z}}\right]^{e^{2y}}{H}(e^z)dz\nonumber\\
&<&4\pi e^{5y}\sum\limits_{n=1}^{\infty}n^4\int_1^{\infty}w^5\left[e^{-\pi n^2w^2}\right]^{e^{2y}}{H}(w)dw\nonumber\\
&&+2\pi \left(\int_{1}^{\infty}{H}(w)dw\right)e^{\frac{9}{2}y}\sum\limits_{n=1}^{\infty}n^4e^{-\pi n^2e^{y}}+2\pi e^{3y}\sum\limits_{n=1}^{\infty}n^4e^{-\pi n^2}\sum_{m=1}^{\infty}8\pi^2\int_{\frac{y}{2}}^{y}(me^z)^4e^{-\pi(me^z)^2}dz\nonumber\\
&<&4\pi \left(\int_1^{\infty}w^5{H}(w)dw\right)e^{5y}\sum\limits_{n=1}^{\infty}n^4e^{-\pi n^2e^{2y}}\nonumber\\
&&+2\pi \left(\int_{1}^{\infty}{H}(w)dw\right)e^{\frac{9}{2}y}\sum\limits_{n=1}^{\infty}n^4e^{-\pi n^2e^{y}}+2\pi e^{3y-2e^y}\sum\limits_{n=1}^{\infty}n^4e^{-\pi n^2}\sum_{m=1}^{\infty}8\pi^2\int_{\frac{y}{2}}^{y}\frac{6}{m^2e^{2z}}dz
\end{eqnarray*}
\begin{eqnarray*}
&<&4\pi \left(\int_1^{\infty}w^5{H}(w)dw\right)e^{5y- 2e^{2y}}\sum\limits_{n=1}^{\infty}n^4e^{-n^2}\nonumber\\
&&+2\pi \left(\int_{1}^{\infty}{H}(w)dw\right)e^{\frac{9}{2}y-2e^y}\sum\limits_{n=1}^{\infty}n^4e^{-n^2}+8\pi^5 y e^{2y-2e^y}\sum\limits_{n=1}^{\infty}n^4e^{-n^2}\nonumber\\
&<&\pi^3 \left(6\int_1^{\infty}w^5{H}(w)dw+8\pi^4\right)e^{5y-2e^{y}}\nonumber\\
&<&\pi^3 \left(\sum_{n=1}^{\infty}\int_1^{\infty}\frac{48\pi^2\cdot 6! (nw)^9}{\pi^6(nw)^{12}}dw+8\pi^4\right)e^{5y-2e^{y}}\nonumber\\
&<&\pi^3 \left(\frac{4\cdot6!}{\pi^2}+8\pi^4\right)e^{5y-2e^{y}}\nonumber\\
&<&24\pi^7e^{5y-2e^{y}}.
\end{eqnarray*}
Then, for $\sigma\in(\frac{1}{2},1)$,
\begin{equation}\label{rem222}
{\cal U}_{\sigma}(y)<96\pi^8e^{5y-2e^{y}},\ \ \ \ y\ge 0.
\end{equation}
\end{rem}

For  $\sigma\in\mathbb{R}$, define
\begin{eqnarray*}
{\cal W}_{\sigma}(x)=\int_{-\infty}^{\infty}R(e^y)R(e^{x-y})e^{x+(2\sigma-1) y}dy,\ \ \ \ x\in\mathbb{R}.
\end{eqnarray*}
We have
$$
{\cal W}_{\sigma}(x)=({\cal R}_{1-\sigma}\ast{\cal R}_{\sigma})(x)e^{\sigma x},\ \ \ \ x\ge0,
$$
where ``$\ast$" is the convolution operator and
$$
{\cal R}_{\sigma}(x):=R(e^{x})e^{\sigma x},\ \ \ \ x\in\mathbb{R}.
$$
It is well known that the following inversion relations hold (cf.  \cite[Page 209]{E}):
\begin{eqnarray}\label{for1}
G(y)=y^{-1}G(y^{-1}),\ \ \ \ H(y)=y^{-1}H(y^{-1}),\ \ \ \ y>0.
\end{eqnarray}
For $y\in\mathbb{R}$, denote by $\lceil y\rceil$ the smallest integer greater than or equal to $y$.

\begin{lem}\label{lem22}
(i) For $n\in\mathbb{Z}$,
$$
\lim_{y\downarrow 0}y^{n}[yG(y)-1]=0.
$$

\noindent (ii) For $n\in\mathbb{N}$,
\begin{eqnarray}\label{28A}
C_n:=\sup_{y>0}\left\{y^nR(y)\right\}<\infty.
\end{eqnarray}

\noindent (iii) For $y>0$,
\begin{eqnarray}\label{for2}
y^{-1}G'(y^{-1})=-y^2G'(y)-G(y^{-1}).
\end{eqnarray}

\noindent (iv)  For $n\in\mathbb{Z}$,
$$
\lim_{y\downarrow 0}y^{n}[y^2G'(y)+1]=0.
$$

\noindent (v) For $\sigma\in\mathbb{R}$,
\begin{eqnarray*}
&&\int_{0}^{\infty}\int_1^{\infty}H(y)R(xy)x^{\sigma-1}y^{2\sigma-1}dxdy\le C_{\lceil |\sigma|\rceil+1}\int_{0}^{\infty}y^{2(\sigma-1)-\lceil |\sigma|\rceil}H(y)dy.
\end{eqnarray*}

\noindent (vi)  For $\sigma\in\mathbb{R}$, \begin{eqnarray*}
{\cal W}_{\sigma}(x)\le (C_1+C_{4\lceil |\sigma|\rceil+3})C_{2\lceil |\sigma|\rceil+3}e^{-2(\lceil |\sigma|\rceil+1)x},\ \ \ \ \forall x\ge0.
\end{eqnarray*}

\noindent (vii) Define
\begin{equation}\label{C}
{\cal C}=\sup_{y>0}\left\{(y^3+y)R(y)\right\}.
\end{equation}
We  have $1\le{\cal C}<\infty$ and for  $\sigma\in(\frac{1}{2},1)$,
$$
{\cal W}_{\sigma}(x)<2{\cal C}^2,\ \ \ \ \forall x\ge0.
$$

\end{lem}

\noindent {\bf Proof.}\ \ (i) By (\ref{for1}), we get
$$
\lim_{y\downarrow 0}y^{n}[yG(y)-1]=\lim_{y\downarrow 0}y^{n}[G(y^{-1})-1]=\lim_{y\downarrow 0}y^{n}R(y^{-1})=0.
$$

\noindent (ii) This is a direct consequence of (i).

\noindent (iii) This is a direct consequence of (\ref{for1}).

\noindent (iv) By (\ref{for2}), we get
\begin{eqnarray*}
\lim_{y\downarrow 0}y^{n}[y^2G'(y)+1]&=&-\lim_{y\downarrow 0}y^{n}[y^{-1}G'(y^{-1})+G(y^{-1})-1]\\
&=&-\lim_{y\downarrow 0}y^{n-1}G'(y^{-1})-\lim_{y\downarrow 0}y^{n}R(y^{-1})\\
&=&0.
\end{eqnarray*}

\noindent (v) We have
\begin{eqnarray*}
\int_{0}^{\infty}\int_1^{\infty}H(y)R(xy)x^{\sigma-1}y^{2\sigma-1}dxdy
&\le&C_{\lceil |\sigma|\rceil+1}\int_{0}^{\infty}\int_1^{\infty}H(y)(xy)^{-(\lceil |\sigma|\rceil+1)}x^{\sigma-1}y^{2\sigma-1}dxdy\\
&\le&C_{\lceil |\sigma|\rceil+1}\int_{0}^{\infty}y^{2(\sigma-1)-\lceil |\sigma|\rceil}H(y)dy.
\end{eqnarray*}

\noindent (vi) For $x\ge0$, we have
\begin{eqnarray*}
{\cal W}_{\sigma}(x)&\le&C_{4\lceil |\sigma|\rceil+3}C_{2\lceil |\sigma|\rceil+3}\int_0^{\infty}e^{-(4\lceil |\sigma|\rceil+3)y}e^{(2\lceil |\sigma|\rceil+3)(y-x)}e^{x+(2\sigma-1) y}dy\\
&&+\int_{0}^1R(z)R(e^xz^{-1})e^xz^{2(\sigma-1)}dz\\
&\le&C_{4\lceil |\sigma|\rceil+3}C_{2\lceil |\sigma|\rceil+3}e^{-2(\lceil |\sigma|\rceil+1)x}\int_{0}^{\infty}e^{-y}dy\\
&&+C_{2\lceil |\sigma|\rceil+3}e^{-2(\lceil |\sigma|\rceil+1)x}\int_{0}^1z^{2(\lceil |\sigma|\rceil+\sigma)+1}R(z)dz\\
&\le&(C_1+C_{4\lceil |\sigma|\rceil+3})C_{2\lceil |\sigma|\rceil+3}e^{-2(\lceil |\sigma|\rceil+1)x}.
\end{eqnarray*}

\noindent (vii) By (i) and (ii), we find that $1\le {\cal C}<\infty$. For $x\ge0$, we have
\begin{eqnarray*}
{\cal W}_{\sigma}(x)&\le&{\cal C}^2\int_0^{\infty}e^{-3y}e^{y-x}e^{x+(2\sigma-1) y}dy+\int_{0}^1R(z)R(e^xz^{-1})e^xz^{2(\sigma-1)}dz\\
&<&{\cal C}^2\int_{0}^{\infty}e^{-y}dy+{\cal C}e^{-2x}\int_{0}^1z^2R(z)dz\\
&<&2{\cal C}^2.
\end{eqnarray*}
Therefore, the proof is complete.\hfill \fbox

\section{Modulus representation of  $\xi$ and proof of Theorem \ref{thm0000}}\setcounter{equation}{0}

For  $\sigma\in\mathbb{R}$, we define
\begin{eqnarray}\label{Jan27a1}
{\cal S}_{\sigma}&=&-[\sigma-2(1-\sigma)^2]R(1)-\sigma R'(1)+\frac{(1-\sigma)(1+2\sigma)R^2(1)}{2}-\sigma R(1)R'(1)\nonumber\\
&&+2\sigma(1-\sigma)(2\sigma-1)\int_{1}^{\infty}y^{-2\sigma}R(y)dy-2\sigma(1-\sigma)(2\sigma-1)\int_{1}^{\infty}y^{1-2\sigma}R(y)dy\nonumber\\
&&+2(1-\sigma)^2(2\sigma-1)\int_{1}^{\infty}y^{2(\sigma-1)}R(y)dy-\sigma(1-\sigma)(2\sigma-1)\int_{1}^{\infty}y^{1-2\sigma}R^2(y)dy\nonumber\\
&&+\frac{2\sigma+1}{2}\int_{1}^{\infty}y^{2\sigma}R^2(y)dy+(1-\sigma)\int_{1}^{\infty}y^{3-2\sigma}[R'(y)]^2dy\nonumber\\
&&-(1-\sigma)\int_{1}^{\infty}y^{1-2\sigma}[yR'(y)]^2dy+\sigma\int_{1}^{\infty}y^{2\sigma-1}[yR'(y)]^2dy\nonumber\\
&&-\sigma(1-\sigma)\int_{0}^{\infty}y^{2(\sigma-1)}R(y) \left\{(1-\sigma)[yR(y)+y-1]+[y^2R'(y)+1]\right\}dy,
\end{eqnarray}
\begin{eqnarray}\label{Jan27a2}
{\cal T}_{\sigma}=-\int_{0}^{\infty}y^{2(\sigma-1)}R(y)\left\{\sigma [yR(y)+y-1]+[y^2R'(y)+1]\right\}dy,
\end{eqnarray}
\begin{eqnarray*}
a_{\sigma}(k)=2^{2k+1}\int_1^{\infty}\int_1^{\infty}(\ln x)^{2k}{\cal J}_{\sigma}(x^2y){\eta}_{\sigma}(y)dxdy,\ \ \ \ k\in\{0\}\cup\mathbb{N},
\end{eqnarray*}
\begin{eqnarray*}
c_{\sigma}(0)&=&\int_1^{\infty}\left[-y^3{\cal J}'''_{\sigma}(y)-\frac{9y^2}{2}{\cal J}''_{\sigma}(y)-\frac{(11-8\sigma^2)y}{4}{\cal J}'_{\sigma}(y)+\frac{1+8\sigma^2}{8}{\cal J}_{\sigma}(y)\right]{\eta}_{\sigma}(y)dy\\
&&+\left(\sigma^2-\frac{1}{4}\right)^2a_{\sigma}(0),
\end{eqnarray*}
\begin{eqnarray*}
c_{\sigma}(1)=\int_1^{\infty}\left[y{\cal J}'_{\sigma}(y)+\frac{1}{2}{\cal J}_{\sigma}(y)\right]{\eta}_{\sigma}(y)dy-\frac{1}{2}\left(\sigma^2-\frac{1}{4}\right)^2a_{\sigma}(1)+2\left(\sigma^2+\frac{1}{4}\right)a_{\sigma}(0),
\end{eqnarray*}
and
\begin{eqnarray*}
c_{\sigma}(k)&=&\frac{(-1)^k}{(2k)!}\Bigg[\left(\sigma^2-\frac{1}{4}\right)^2a_{\sigma}(k)-4k(2k-1)\left(\sigma^2+\frac{1}{4}\right)a_{\sigma}(k-1)\\
&&\quad\quad\quad\quad+2k(2k-1)(2k-2)(2k-3)a_{\sigma}(k-2) \Bigg],\ \ \ \ k\ge 2.
\end{eqnarray*}

In this section, we will prove the following representation of the modulus of the function $\xi$. 

\begin{thm}\label{thm000} (i) For $\sigma\in\mathbb{R}$ and $t\in\mathbb{R}$,
\begin{eqnarray}\label{13xyz}
2|\xi(\sigma-it)|^2={\cal S}_{\sigma}+{\cal T}_{\sigma}t^2+\{t^2+(1-\sigma)^2\}\{t^2+\sigma^2\}\int_0^{\infty}{\cal W}_{\sigma}(x)e^{-\sigma x}\cos(tx)dx.
\end{eqnarray}

\noindent (ii)  For $\sigma\in\mathbb{R}$,
\begin{eqnarray*}
{\cal S}_{\sigma}&=&\int_1^{\infty}\Bigg\{-2y^3{\cal J}'''_{\sigma-\frac{1}{2}}(y)-9y^2{\cal J}''_{\sigma-\frac{1}{2}}(y)+y\left[4\left(\sigma-\frac{1}{2}\right)^2-\frac{11}{2}\right]{\cal J}'_{\sigma-\frac{1}{2}}(y)\\
&&\quad\quad\quad+\left[2\left(\sigma-\frac{1}{2}\right)^2+\frac{1}{4}\right]{\cal J}_{\sigma-\frac{1}{2}}(y)\Bigg\}{\eta}_{\sigma-\frac{1}{2}}(y)dy,
\end{eqnarray*}
\begin{eqnarray*}
{\cal T}_{\sigma}=\int_1^{\infty}\left[2y{\cal J}'_{\sigma-\frac{1}{2}}(y)+{\cal J}_{\sigma-\frac{1}{2}}(y)\right]{\eta}_{\sigma-\frac{1}{2}}(y)dy,
\end{eqnarray*}
and
\begin{eqnarray}\label{Mar28as}
\int_0^{\infty}{\cal W}_{\sigma}(x)e^{-\sigma x}\cos(tx)dx=4\int_1^{\infty}\int_1^{\infty}\cos(2t\ln x){\cal J}_{\sigma-\frac{1}{2}}(x^2y){\eta}_{\sigma-\frac{1}{2}}(y)dxdy.
\end{eqnarray}

\noindent (iii) For $\sigma\in\mathbb{R}$ and $t\in\mathbb{R}$,
\begin{equation}\label{Feb11}
\left|\xi\left(\sigma-it\right)\right|^2=\sum_{k=0}^{\infty}c_{\sigma-\frac{1}{2}}(k)t^{2k}.
\end{equation}
\end{thm}

It is well known that $\xi(\sigma-it)\not=0$ for any $t\in\mathbb{R}$ if $|\sigma-\frac{1}{2}|\ge\frac{1}{2}$. Note that 
$$\{t^2+(1-\sigma)^2\}\{t^2+\sigma^2\}=\left[t^2+\left(\sigma-\frac{1}{2}\right)^2+\frac{1}{4}\right]^2-\left(\sigma-\frac{1}{2}\right)^2.
$$
Let $\tau=\sigma-\frac{1}{2}$. Then, Theorem \ref{thm0000} is a direct consequence of Theorem \ref{thm000} (i) and (ii). 
\vskip 0.5cm
\noindent {\bf Proof of Theorem \ref{thm000} (i).}\ \ By (\ref{Oct22a}) and (\ref{for1}), we find that for $x\ge0$,
\begin{eqnarray*}
{\cal U}_{\sigma}(x)&=&e^{(1-\sigma)x}\int_{0}^{\infty}H(e^{x+y})H(e^y)e^{2(1-\sigma)y}dy+e^{(1-\sigma)x}\int_{0}^{x}H(e^{x-y})H(e^y)e^{(2\sigma-1)y}dy\nonumber\\
&&+e^{\sigma x}\int_{0}^{\infty}H(e^{x+y})H(e^y)e^{2\sigma y}dy\\
&=&e^{(1-\sigma)x}\int_{1}^{\infty}H(y)H(e^xy)y^{1-2\sigma}dy+e^{(1-\sigma)x}\int_{e^{-x}}^{1}H(y^{-1})H(e^xy)y^{-2\sigma}dy\nonumber\\
&&+e^{\sigma x}\int_{1}^{\infty}H(y)H(e^xy)y^{2\sigma-1}dy\nonumber\\
&=&e^{(1-\sigma)x}\int_{1}^{\infty}H(y)H(e^xy)y^{1-2\sigma}dy+e^{(1-\sigma)x}\int_{e^{-x}}^{1}H(y)H(e^xy)y^{1-2\sigma}dy\nonumber\\
&&+e^{\sigma x}\int_{1}^{\infty}H(y)H(e^xy)y^{2\sigma-1}dy\nonumber\\
&=&e^{\sigma x}\int_{1}^{\infty}H(y)H(e^xy)y^{2\sigma-1}dy+e^{(1-\sigma)x}\int_{e^{-x}}^{\infty}H(y)H(e^xy)y^{1-2\sigma}dy\nonumber\\
&=&e^{\sigma x}\int_{1}^{\infty}H(y)H(e^xy)y^{2\sigma-1}dy+e^{(\sigma-1)x}\int_{1}^{\infty}H(y)H(e^{-x}y)y^{1-2\sigma}dy.
\end{eqnarray*}
Then, by (\ref{RRR2}), we get
\begin{eqnarray}\label{12/12d1}
|\xi(\sigma-it)|^2
&=&\frac{1}{2}\int_0^{\infty}\cos(tx)\Bigg[e^{\sigma x}\int_{1}^{\infty}H(y)H(e^xy)y^{2\sigma-1}dy\nonumber\\
&&\quad\quad\quad\quad+e^{(\sigma-1)x}\int_{1}^{\infty}H(y)H(e^{-x}y)y^{1-2\sigma}dy\Bigg]dx.
\end{eqnarray}

By integration by parts, (\ref{for1}), (\ref{for2}) and Lemma \ref{lem22} (v), we obtain that
\begin{eqnarray*}
&&\int_0^{\infty}\cos(tx)e^{\sigma x}\int_{1}^{\infty}H(y)H(e^xy)y^{2\sigma-1}dydx\\
&=&\int_{1}^{\infty}y^{2\sigma-1}H(y)\int_0^{\infty}\cos(tx)e^{\sigma x}H(e^xy)dxdy\\
&=&\int_{1}^{\infty}y^{\sigma-1}H(y)\int_y^{\infty}\frac{\left(\frac{z}{y}\right)^{it}+\left(\frac{z}{y}\right)^{-it}}{2z^{1-\sigma}}H(z)dzdy\\
&=&\int_{1}^{\infty}y^{\sigma-1}H(y)\int_y^{\infty}\frac{z^{\sigma}\left[\left(\frac{z}{y}\right)^{it}+\left(\frac{z}{y}\right)^{-it}\right]}{2}d[zR(z)]'dy
\end{eqnarray*}
\begin{eqnarray*}
&=&-\int_{1}^{\infty}y^{2\sigma-1}H(y)\{R(y)+yR'(y)\}dy\\
&&-\int_{1}^{\infty}y^{\sigma-1}H(y)\int_y^{\infty}\frac{z^{\sigma-1}\left[\left(\frac{z}{y}\right)^{it}(\sigma+it)+\left(\frac{z}{y}\right)^{-it}(\sigma-it)\right]}{2}d[zR(z)]dy\\
&=&\int_{1}^{\infty}y^{2\sigma-1}H(y)\{(\sigma-1)R(y)-yR'(y)\}dy+\int_{1}^{\infty}y^{\sigma-1}H(y)\int_y^{\infty}\\
&&\ \ \ \ \frac{z^{\sigma-1}R(z)\left[\left(\frac{z}{y}\right)^{it}(\sigma+it)(\sigma-1+it)+\left(\frac{z}{y}\right)^{-it}(\sigma-it)(\sigma-1-it)\right]}{2}dzdy\\
&=&\int_{1}^{\infty}y^{2\sigma-1}H(y)\{(\sigma-1)R(y)-yR'(y)\}dy+\int_{1}^{\infty}y^{2\sigma-1}H(y)\int_1^{\infty}\\
&&\ \ \ \ \frac{R(xy)\left[x^{it}(\sigma+it)(\sigma-1+it)+x^{-it}(\sigma-it)(\sigma-1-it)\right]}{2x^{1-\sigma}}dxdy,
\end{eqnarray*}
and
\begin{eqnarray*}
&&\int_0^{\infty}\frac{\cos(tx)}{e^{(1-\sigma)x}}\int_{1}^{\infty}H(y)H(e^{-x}y)y^{1-2\sigma}dydx\\
&=&\int_{1}^{\infty}y^{1-2\sigma}H(y)\int_0^{\infty}\frac{\cos(tx)}{e^{(1-\sigma)x}}H(e^{-x}y)dxdy\\
&=&\int_{1}^{\infty}y^{-\sigma}H(y)\int_0^y\frac{z^{-\sigma}\left[\left(\frac{z}{y}\right)^{it}+\left(\frac{z}{y}\right)^{-it}\right]}{2}H(z)dzdy\\
&=&\int_{1}^{\infty}y^{-\sigma}H(y)\int_{\frac{1}{y}}^{\infty}\frac{x^{\sigma-1}\left[(xy)^{it}+(xy)^{-it}\right]}{2}H(x)dxdy\\
&=&\int_{1}^{\infty}y^{-\sigma}H(y)\int_{\frac{1}{y}}^{\infty}\frac{x^{\sigma}\left[(xy)^{it}+(xy)^{-it}\right]}{2}d[xR(x)]'dy\\
&=&-\int_{1}^{\infty}y^{-2\sigma}H(y)\{R(y^{-1})+y^{-1}R'(y^{-1})\}dy\\
&&-\int_{1}^{\infty}y^{-\sigma}H(y)\int_{\frac{1}{y}}^{\infty}\frac{x^{\sigma-1}\left[(xy)^{it}(\sigma+it)+(xy)^{-it}(\sigma-it)\right]}{2}d[xR(x)]dy\\
&=&\int_{1}^{\infty}y^{1-2\sigma}H(y)\{\sigma [R(y)+1]+yR'(y)+(1-\sigma)y^{-1}\}dy+\int_{1}^{\infty}y^{-\sigma}H(y)\int_{\frac{1}{y}}^{\infty}\\
&&\ \ \ \ \frac{x^{\sigma-1}R(x)[(xy)^{it}(\sigma-1+it)(\sigma+it)+(xy)^{-it}(\sigma-1-it)(\sigma-it)]}{2}dxdy\\
&=&\int_{1}^{\infty}y^{1-2\sigma}H(y)\{\sigma [R(y)+1]+yR'(y)+(1-\sigma)y^{-1}\}dy+\int_{1}^{\infty}y^{-2\sigma}H(y)\int_1^{\infty}\\
&&\ \ \ \ \frac{R\left(\frac{z}{y}\right)[z^{it}(\sigma-1+it)(\sigma+it)+z^{-it}(\sigma-1-it)(\sigma-it)]}{2z^{1-\sigma}}dzdy
\end{eqnarray*}
\begin{eqnarray*}
&=&\int_{1}^{\infty}y^{1-2\sigma}H(y)\{\sigma [R(y)+1]+yR'(y)+(1-\sigma)y^{-1}\}dy+\int_{0}^{1}y^{2\sigma-1}H(y)\int_1^{\infty}\\
&&\ \ \ \ \frac{R(xy)\left[x^{it}(\sigma-1+it)(\sigma+it)+x^{-it}(\sigma-1-it)(\sigma-it)\right]}{2x^{1-\sigma}}dxdy.
\end{eqnarray*}
Then, by (\ref{12/12d1}), we get
\begin{eqnarray}\label{13a}
&&2|\xi(\sigma-it)|^2\nonumber\\
&=&\int_{1}^{\infty}y^{2\sigma-1}H(y)\{(\sigma-1)R(y)-yR'(y)\}dy\nonumber\\
&&+\int_{1}^{\infty}y^{1-2\sigma}H(y)\{\sigma [R(y)+1]+yR'(y)+(1-\sigma)y^{-1}\}dy\\
&&+\int_{0}^{\infty}y^{2\sigma-1}H(y)\int_1^{\infty}\frac{R(xy)\left[x^{it}(\sigma-1+it)(\sigma+it)+x^{-it}(\sigma-1-it)(\sigma-it)\right]}{2x^{1-\sigma}}dxdy.\nonumber
\end{eqnarray}

By Lemma \ref{lem22}, we get
\begin{eqnarray*}
&&\int_{0}^{\infty}\int_1^{\infty}\frac{y^{2\sigma-1}H(y)R(xy)\left[x^{it}(\sigma-1+it)(\sigma+it)+x^{-it}(\sigma-1-it)(\sigma-it)\right]}{2x^{1-\sigma}}dxdy\\
&=&\int_{0}^{\infty}z^{\sigma-1}R(z)\int_0^{z}\frac{y^{\sigma-1}H(y)\left[\left(\frac{z}{y}\right)^{it}(\sigma-1+it)(\sigma+it)+\left(\frac{z}{y}\right)^{-it}(\sigma-1-it)(\sigma-it)\right]}{2}dydz\\
&=&\int_{0}^{\infty}z^{\sigma-1}R(z)\int_{\frac{1}{z}}^{\infty}\frac{y^{-\sigma}H(y)\left[(yz)^{it}(\sigma-1+it)(\sigma+it)+(yz)^{-it}(\sigma-1-it)(\sigma-it)\right]}{2}dydz\\
&=&\int_{0}^{\infty}z^{\sigma-1}R(z)\int_{\frac{1}{z}}^{\infty}\frac{y^{1-\sigma}\left[(yz)^{it}(\sigma-1+it)(\sigma+it)+(yz)^{-it}(\sigma-1-it)(\sigma-it)\right]}{2}d[yR(y)]'dz\\
&=&-[(\sigma-1)\sigma-t^2]\int_{0}^{\infty}z^{2(\sigma-1)}R(z)[R(z^{-1})+z^{-1}R'(z^{-1})]dz+\int_{0}^{\infty}z^{\sigma-1}R(z)\int_{\frac{1}{z}}^{\infty}\\
&&\ \ \ \ \frac{y^{-\sigma}\left[(yz)^{it}\{(\sigma-1)^2+t^2\}(\sigma+it)+(yz)^{-it}\{(\sigma-1)^2+t^2\}(\sigma-it)\right]}{2}d[yR(y)]dz\\
&=&\int_{0}^{\infty}z^{2(\sigma-1)}R(z)\left\{[\sigma(\sigma-1)-t^2][z^2R'(z)+1]-\sigma[(\sigma-1)^2+t^2][zR(z)+z-1]\right\}dz\\
&&+\int_{0}^{\infty}z^{\sigma-1}R(z)\int_{\frac{1}{z}}^{\infty}\frac{y^{-\sigma}R(y)\left[(yz)^{it}\{(\sigma-1)^2+t^2\}\{\sigma^2+t^2\}+(yz)^{-it}\{(\sigma-1)^2+t^2\}\{\sigma^2+t^2\}\right]}{2}dydz\\
&=&\int_{0}^{\infty}z^{2(\sigma-1)}R(z)\Bigg\{\sigma(\sigma-1) \{[z^2R'(z)+1]+(1-\sigma)[zR(z)+z-1]\}\\
&&\ \ \ \ -\{\sigma [zR(z)+z-1]+[z^2R'(z)+1]\}t^2\Bigg\}dz\\
&&+\int_{0}^{\infty}\int_1^{\infty}\frac{z^{2(\sigma-1)}R(z)R\left(\frac{x}{z}\right)[x^{it}\{(\sigma-1)^2+t^2\}\{\sigma^2+t^2\}+x^{-it}\{(\sigma-1)^2+t^2\}\{\sigma^2+t^2\}]}{2x^{\sigma}}dxdz
\end{eqnarray*}
\begin{eqnarray*}
&=&\int_{0}^{\infty}y^{2(\sigma-1)}R(y)\Bigg\{-\sigma(1-\sigma) \{(1-\sigma)[yR(y)+y-1]+[y^2R'(y)+1]\}\\
&&\ \ \ \ -t^2\{\sigma [yR(y)+y-1]+[y^2R'(y)+1]\}\Bigg\}dy\\
&&+\{(\sigma-1)^2+t^2\}\{\sigma^2+t^2\}\int_{0}^{\infty}\int_1^{\infty}\frac{y^{-2\sigma}R(y^{-1})R(xy)\cos(t\ln x)}{x^{\sigma}}dxdy\\
&=&\int_{0}^{\infty}y^{2(\sigma-1)}R(y)\Bigg\{-\sigma(1-\sigma) \{(1-\sigma)[yR(y)+y-1]+[y^2R'(y)+1]\}\\
&&\ \ \ \ -t^2\{\sigma [yR(y)+y-1]+[y^2R'(y)+1]\}\Bigg\}dy\\
&&+\{t^2+(1-\sigma)^2\}\{t^2+\sigma^2\}\int_{-\infty}^{\infty}\int_0^{\infty}R(e^y)R(e^{x-y})e^{(1-\sigma) x}e^{(2\sigma-1) y}\cos(tx)dxdy.
\end{eqnarray*}
Then, by (\ref{13a}), we obtain that
\begin{eqnarray}\label{diff}
2|\xi(\sigma-it)|^2
&=&\int_{1}^{\infty}y^{2\sigma-1}H(y)\{(\sigma-1)R(y)-yR'(y)\}dy\nonumber\\
&&+\int_{1}^{\infty}y^{1-2\sigma}H(y)\{\sigma [R(y)+1]+yR'(y)+(1-\sigma)y^{-1}\}dy\nonumber\\
&&-\sigma(1-\sigma)\int_{0}^{\infty}y^{2(\sigma-1)}R(y) \{(1-\sigma)[yR(y)+y-1]+[y^2R'(y)+1]\}dy\nonumber\\
&&+{\cal T}_{\sigma}t^2+\{t^2+(1-\sigma)^2\}\{t^2+\sigma^2\}\int_0^{\infty}{\cal W}_{\sigma}(x)e^{-\sigma x}\cos(tx)dx.
\end{eqnarray}

We have
\begin{eqnarray*}
&&\int_{1}^{\infty}y^{2\sigma-1}H(y)\{(\sigma-1)R(y)-yR'(y)\}dy\\
&=&(\sigma-1)\int_{1}^{\infty}y^{2\sigma}R(y)d[yR(y)]'-\int_{1}^{\infty}y^{2\sigma+1}R'(y)d[yR(y)]'\\
&=&(\sigma-1)\int_{1}^{\infty}y^{2\sigma-1}d[R(y)+yR'(y)]-\int_{1}^{\infty}y^{2\sigma+1}R'(y)d[R(y)+yR'(y)]\\
&=&(1-\sigma)[R(1)+R'(1)]+(1-\sigma)(2\sigma-1)\int_{1}^{\infty}[R(y)+yR'(y)]y^{2(\sigma-1)}dy\\
&&-\frac{1}{2}\int_{1}^{\infty}y^{2\sigma+1}dR^2(y)-\frac{1}{2}\int_{1}^{\infty}y^{2\sigma}d[yR'(y)]^2\\
&=&(1-\sigma)[R(1)+R'(1)]+(1-\sigma)(2\sigma-1)\int_{1}^{\infty}y^{2(\sigma-1)}R(y)dy\\
&&+(1-\sigma)(2\sigma-1)\int_{1}^{\infty}y^{2\sigma-1}dR(y)\\
&&+\frac{R^2(1)}{2}+\frac{2\sigma+1}{2}\int_{1}^{\infty}y^{2\sigma}R^2(y)dy+\frac{[R'(1)]^2}{2}+\sigma\int_{1}^{\infty}y^{2\sigma-1}[yR'(y)]^2dy
\end{eqnarray*}
\begin{eqnarray*}
&=&(1-\sigma)[R(1)+R'(1)]+(1-\sigma)(2\sigma-1)\int_{1}^{\infty}y^{2(\sigma-1)}R(y)dy\\
&&-(1-\sigma)(2\sigma-1)\left[R(1)+(2\sigma-1)\int_{1}^{\infty}y^{2(\sigma-1)}R(y)dy\right]\\
&&+\frac{R^2(1)}{2}+\frac{2\sigma+1}{2}\int_{1}^{\infty}y^{2\sigma}R^2(y)dy+\frac{[R'(1)]^2}{2}+\sigma\int_{1}^{\infty}y^{2\sigma-1}[yR'(y)]^2dy\\
&=&2(1-\sigma)^2R(1)+(1-\sigma)R'(1)+\frac{R^2(1)}{2}+\frac{[R'(1)]^2}{2}\\
&&+2(1-\sigma)^2(2\sigma-1)\int_{1}^{\infty}y^{2(\sigma-1)}R(y)dy+\frac{2\sigma+1}{2}\int_{1}^{\infty}y^{2\sigma}R^2(y)dy+\sigma\int_{1}^{\infty}y^{2\sigma-1}[yR'(y)]^2dy,
\end{eqnarray*}
and
\begin{eqnarray*}
&&\int_{1}^{\infty}y^{1-2\sigma}H(y)\{\sigma [R(y)+1]+yR'(y)+(1-\sigma)y^{-1}\}dy\\
&=&\sigma\int_{1}^{\infty}y^{2(1-\sigma)}R(y)d[R(y)+yR'(y)]+\sigma\int_{1}^{\infty}y^{2(1-\sigma)}d[R(y)+yR'(y)]\\
&&+\int_{1}^{\infty}y^{3-2\sigma}R'(y)d[R(y)+yR'(y)]+(1-\sigma)\int_{1}^{\infty}y^{1-2\sigma}d[R(y)+yR'(y)]\\
&=&\frac{\sigma}{2}\int_{1}^{\infty}y^{2(1-\sigma)}dR^2(y)-\sigma R(1)R'(1)-\sigma\int_{1}^{\infty}yR'(y)d[y^{2(1-\sigma)}R(y)]\\
&&-\sigma[R(1)+R'(1)]-2\sigma(1-\sigma)\int_{1}^{\infty}[R(y)+yR'(y)]y^{1-2\sigma}dy\\
&&+\int_{1}^{\infty}y^{3-2\sigma}[R'(y)]^2dy+\frac{1}{2}\int_{1}^{\infty}y^{2(1-\sigma)}d[yR'(y)]^2\\
&&-(1-\sigma)\left[R(1)+(1-2\sigma)\int_{1}^{\infty}y^{-2\sigma}R(y)dy\right]\\
&&-(1-\sigma)\left[R'(1)+(1-2\sigma)\int_{1}^{\infty}y^{1-2\sigma}dR(y)\right]\\
&=&-\frac{\sigma R^2(1)}{2}-\sigma(1-\sigma)\int_{1}^{\infty}y^{1-2\sigma}R^2(y)dy\\
&&-\sigma R(1)R'(1)+\sigma(1-\sigma)R^2(1)+2\sigma(1-\sigma)^2\int_{1}^{\infty}y^{1-2\sigma}R^2(y)dy\\
&&-\sigma\int_{1}^{\infty}y^{3-2\sigma}[R'(y)]^2dy-\sigma[R(1)+R'(1)]-2\sigma(1-\sigma)\int_{1}^{\infty}y^{1-2\sigma}R(y)dy\\
&&+2\sigma(1-\sigma)R(1)+4\sigma(1-\sigma)^2\int_{1}^{\infty}y^{1-2\sigma}R(y)dy\\
&&+\int_{1}^{\infty}y^{3-2\sigma}[R'(y)]^2dy-\frac{[R'(1)]^2}{2}-(1-\sigma)\int_{1}^{\infty}y^{1-2\sigma}[yR'(y)]^2dy\\
&&-(1-\sigma)\left[R(1)+(1-2\sigma)\int_{1}^{\infty}y^{-2\sigma}R(y)dy\right]\\
&&-(1-\sigma)\left[R'(1)-(1-2\sigma)R(1)-(1-2\sigma)^2\int_{1}^{\infty}y^{-2\sigma}R(y)dy\right]\end{eqnarray*}
\begin{eqnarray*}
&=&-\sigma R(1)-R'(1)-\frac{\sigma(2\sigma-1)R^2(1)}{2}-\sigma R(1)R'(1)-\frac{[R'(1)]^2}{2}\\
&&+2\sigma(1-\sigma)(2\sigma-1)\int_{1}^{\infty}y^{-2\sigma}R(y)dy-2\sigma(1-\sigma)(2\sigma-1)\int_{1}^{\infty}y^{1-2\sigma}R(y)dy\\
&&-\sigma(1-\sigma)(2\sigma-1)\int_{1}^{\infty}y^{1-2\sigma}R^2(y)dy\\
&&+(1-\sigma)\int_{1}^{\infty}y^{3-2\sigma}[R'(y)]^2dy-(1-\sigma)\int_{1}^{\infty}y^{1-2\sigma}[yR'(y)]^2dy.
\end{eqnarray*}
Then, by summing them up, we get
\begin{eqnarray}\label{diff2}
&&\int_{1}^{\infty}y^{2\sigma-1}H(y)\{(\sigma-1)R(y)-yR'(y)\}dy\nonumber\\
&&+\int_{1}^{\infty}y^{1-2\sigma}H(y)\{\sigma [R(y)+1]+yR'(y)+(1-\sigma)y^{-1}\}dy\nonumber\\
&=&-[\sigma-2(1-\sigma)^2]R(1)-\sigma R'(1)+\frac{(1-\sigma)(1+2\sigma)R^2(1)}{2}-\sigma R(1)R'(1)\nonumber\\
&&+2\sigma(1-\sigma)(2\sigma-1)\int_{1}^{\infty}y^{-2\sigma}R(y)dy-2\sigma(1-\sigma)(2\sigma-1)\int_{1}^{\infty}y^{1-2\sigma}R(y)dy\nonumber\\
&&+2(1-\sigma)^2(2\sigma-1)\int_{1}^{\infty}y^{2(\sigma-1)}R(y)dy-\sigma(1-\sigma)(2\sigma-1)\int_{1}^{\infty}y^{1-2\sigma}R^2(y)dy\nonumber\\
&&+\frac{2\sigma+1}{2}\int_{1}^{\infty}y^{2\sigma}R^2(y)dy+(1-\sigma)\int_{1}^{\infty}y^{3-2\sigma}[R'(y)]^2dy\nonumber\\
&&-(1-\sigma)\int_{1}^{\infty}y^{1-2\sigma}[yR'(y)]^2dy+\sigma\int_{1}^{\infty}y^{2\sigma-1}[yR'(y)]^2dy.
\end{eqnarray}
Thus, (\ref{13xyz}) holds by (\ref{diff}) and (\ref{diff2}).
\vskip 0.5cm
\noindent {\bf Proof of Theorem \ref{thm000} (ii).}\ \ For $\sigma\in\mathbb{R}$, define
\begin{eqnarray}\label{Jan23a}
F_{\sigma}(\lambda)&=&e^{-2\lambda}\int_0^{\infty}\frac{\lambda e^{-\lambda y}}{2\sqrt{y^2+4y}}\Bigg\{\frac{1}{[y+2+\sqrt{y^2+4y}]^{\sigma-\frac{1}{2}}}\nonumber\\
&&\quad\quad\quad\quad\quad\quad\quad\quad\quad\quad+\frac{1}{[y+2-\sqrt{y^2+4y}]^{\sigma-\frac{1}{2}}}\Bigg\}dy,\ \ \ \ \lambda>0,
\end{eqnarray}
\begin{eqnarray}\label{Jan28D1}
{\cal G}_{\sigma}(\lambda)=\sum_{m=1}^{\infty}\sum_{n=1}^{\infty}m^{-\sigma-\frac{1}{2}}n^{\sigma-\frac{3}{2}}F_{\sigma}(\pi mn\lambda),\ \ \ \ \lambda>0,
\end{eqnarray}
and
\begin{eqnarray}\label{Jan28D2}
{\cal H}_{\sigma}(x)={\cal G}_{\sigma}(e^x)e^{-\frac{1}{2}x},\ \ \ \ x\in\mathbb{R}.
\end{eqnarray}
For $x\ge0$, we have that
\begin{eqnarray}\label{13xyzw}
{\cal W}_{\sigma}(x)&=&4e^x\sum_{m=1}^{\infty}\sum_{n=1}^{\infty}\int_{-\infty}^{\infty}\frac{e^{(2\sigma-1) y}}{e^{\pi[(me^y)^2+(ne^{x-y})^2]}}dy\nonumber\\
&=&4e^x\sum_{m=1}^{\infty}\sum_{n=1}^{\infty}\int_0^{\infty}\frac{1}{u^{2\sigma}e^{\pi[(mu^{-1})^2+(ne^xu)^2]}}du\nonumber\\
&=&4\pi^{\sigma-\frac{1}{2}}e^{2\sigma x}\sum_{m=1}^{\infty}\sum_{n=1}^{\infty}n^{2\sigma-1}\int_0^{\infty}\frac{1}{v^{2\sigma}e^{(\pi mne^{x})^2v^{-2}+v^2}}dv\nonumber\\
&=&2\pi^{\sigma-\frac{1}{2}}e^{2\sigma x}\sum_{m=1}^{\infty}\sum_{n=1}^{\infty}n^{2\sigma-1}\int_0^{\infty}\frac{e^{-y^{-1}}}{y^{\frac{3}{2}-\sigma}}\cdot e^{-(\pi mne^{x})^2y}dy\nonumber\\
&=&2\pi^{\sigma-\frac{1}{2}}e^{2\sigma x}\sum_{m=1}^{\infty}\sum_{n=1}^{\infty}n^{2\sigma-1}\Bigg[\int_0^{\left(\pi mne^x\right)^{-1}}+\int_{(\pi mne^x)^{-1}}^{\infty}\Bigg]\frac{1}{y^{\frac{3}{2}-\sigma}}\cdot e^{-(\pi mne^x)\left\{\left(\pi mne^xy\right)^{-1}+\pi mne^xy\right\}}dy\nonumber\\
&=&2\pi^{\sigma-\frac{1}{2}}e^{2\sigma x}\sum_{m=1}^{\infty}\sum_{n=1}^{\infty}n^{2\sigma-1}\int_2^{\infty}\frac{[v-\sqrt{v^2-4}]^{\sigma-\frac{1}{2}}+[v+\sqrt{v^2-4}]^{\sigma-\frac{1}{2}}}{(2\pi mne^x)^{\sigma-\frac{1}{2}}\sqrt{v^2-4}}\cdot e^{-(\pi mne^x)v}dv\nonumber\\
&=&2^{2\sigma}\pi^{\sigma-\frac{1}{2}}e^{2\sigma x}\sum_{m=1}^{\infty}\sum_{n=1}^{\infty}\frac{n^{2\sigma-1}}{(2\pi mne^x)^{\sigma-\frac{1}{2}}}\int_2^{\infty}\frac{e^{-(\pi mne^x)v}}{\sqrt{v^2-4}}\nonumber\\
&&\ \ \ \ \ \ \ \ \ \ \ \ \ \ \ \ \ \ \ \ \cdot\left\{\frac{1}{[v+\sqrt{v^2-4}]^{\sigma-\frac{1}{2}}}+\frac{1}{[v-\sqrt{v^2-4}]^{\sigma-\frac{1}{2}}}\right\}dv\nonumber\\
&=&2^{\sigma+\frac{1}{2}}\pi^{-1}e^{(\sigma-\frac{1}{2})x}\sum_{m=1}^{\infty}\sum_{n=1}^{\infty}\frac{n^{2\sigma-1}}{(mn)^{\sigma+\frac{1}{2}}}\cdot e^{-2\pi mne^x}\int_0^{\infty}\frac{(\pi mne^x)e^{-(\pi mne^x)y}}{\sqrt{y^2+4y}}\nonumber\\
&&\ \ \ \ \ \ \ \ \ \ \ \ \ \ \ \ \ \ \ \ \ \ \ \ \cdot\left\{\frac{1}{[y+2+\sqrt{y^2+4y}]^{\sigma-\frac{1}{2}}}+\frac{1}{[y+2-\sqrt{y^2+4y}]^{\sigma-\frac{1}{2}}}\right\}dy\nonumber\\
&=&2^{\sigma+\frac{3}{2}}\pi^{-1}e^{(\sigma-\frac{1}{2})x}\sum_{m=1}^{\infty}\sum_{n=1}^{\infty}m^{-\sigma-\frac{1}{2}}n^{\sigma-\frac{3}{2}}F_{\sigma}(\pi mne^x)\nonumber\\
&=&2^{\sigma+\frac{3}{2}}\pi^{-1}e^{(\sigma-\frac{1}{2})x}{\cal G}_{\sigma}(e^x)\nonumber\\
&=&2^{\sigma+\frac{3}{2}}\pi^{-1}{\cal H}_{\sigma}(x)e^{\sigma x}.
\end{eqnarray}

By integration by parts, for $t\not=0$, we get
\begin{eqnarray*}
\int_{0}^{\infty}{\cal H}_{\sigma}(x)\cos(tx)dx
&=&\frac{1}{t}\int_{0}^{\infty}{\cal H}_{\sigma}(x)d\sin(tx)\\
&=&\frac{1}{t^2}\int_{0}^{\infty}{\cal H}'_{\sigma}(x)d\cos(tx)\\
&=&-\frac{1}{t^2}{\cal H}'_{\sigma}(0)-\frac{1}{t^2}\int_{0}^{\infty}{\cal H}''_{\sigma}(x)\cos(tx)dx\\
&=&-\frac{1}{t^2}{\cal H}'_{\sigma}(0)-\frac{1}{t^4}\int_{0}^{\infty}{\cal H}'''_{\sigma}(x)d\cos(tx)
\end{eqnarray*}
\begin{eqnarray*}
&=&-\frac{1}{t^2}{\cal H}'_{\sigma}(0)+\frac{1}{t^4}{\cal H}'''_{\sigma}(0)-\frac{1}{t^5}\int_{0}^{\infty}{\cal H}^{(5)}_{\sigma}(x)\sin(tx)dx.
\end{eqnarray*}
Then, by (\ref{13xyz}) and (\ref{13xyzw}), we obtain that  for $t\not=0$,
\begin{eqnarray}\label{Jan1t}
2|\xi(\sigma-it)|^2
&=&{\cal S}_{\sigma}+{\cal T}_{\sigma}t^2+2^{\sigma+\frac{3}{2}}\pi^{-1}\{t^2+(1-\sigma)^2\}\{t^2+\sigma^2\}\nonumber\\
&&\ \ \ \ \cdot\Bigg\{-\frac{1}{t^2}{\cal H}'_{\sigma}(0)+\frac{1}{t^4}{\cal H}_{\sigma}'''(0)-\frac{1}{t^5}\int_{0}^{\infty}{\cal H}^{(5)}_{\sigma}(x)\sin(tx)dx\Bigg\}\nonumber\\
&=&{\cal S}_{\sigma}-2^{\sigma+\frac{3}{2}}\pi^{-1}\left\{[\sigma^2+(1-\sigma)^2]{\cal H}'_{\sigma}(0)-{\cal H}'''_{\sigma}(0)\right\}\nonumber\\
&&\ \ \ \ +\left\{{\cal T}_{\sigma}-2^{\sigma+\frac{3}{2}}\pi^{-1}{\cal H}'_{\sigma}(0)\right\}t^2+\frac{\vartheta_{\sigma}(t)}{t},
\end{eqnarray}
where $\vartheta_{\sigma}(t)$ is a function that is bounded for $|t|\ge 1$.

By Theorem \ref{thm00}, we have
$$
\Xi_{\sigma}(t)=\int_{-\infty}^{\infty}e^{ity}P_{\sigma}(y)dy.
$$
Then, by the Riemann-Lebesgue lemma, we get
$$
\lim_{|t|\rightarrow\infty}\xi(\sigma-it)=0.
$$
Thus, by (\ref{Jan1t}), we obtain that
\begin{eqnarray}\label{Janm1}
{\cal S}_{\sigma}=2^{\sigma+\frac{3}{2}}\pi^{-1}\left\{[\sigma^2+(1-\sigma)^2]{\cal H}'_{\sigma}(0)-{\cal H}'''_{\sigma}(0)\right\},\ \ \ \
{\cal T}_{\sigma}=2^{\sigma+\frac{3}{2}}\pi^{-1}{\cal H}'_{\sigma}(0).
\end{eqnarray}

For $\lambda>0$ and $x\in\mathbb{R}$, by (\ref{Jan23a}) and (\ref{Jan28D2}),  we get
\begin{eqnarray}\label{Jan10a}
&&F'_{\sigma}(\lambda)=e^{-2\lambda}\int_0^{\infty}\frac{(1-2\lambda-\lambda y)e^{-\lambda y}}{2\sqrt{y^2+4y}}\left\{\frac{1}{[y+2+\sqrt{y^2+4y}]^{\sigma-\frac{1}{2}}}+\frac{1}{[y+2-\sqrt{y^2+4y}]^{\sigma-\frac{1}{2}}}\right\}dy,\nonumber\\
&&F''_{\sigma}(\lambda)=e^{-2\lambda}\int_0^{\infty}\frac{-(y+2)(2-2\lambda-\lambda y)e^{-\lambda y}}{2\sqrt{y^2+4y}}\nonumber\\
&&\quad\quad\quad\quad\quad\quad\quad\cdot\left\{\frac{1}{[y+2+\sqrt{y^2+4y}]^{\sigma-\frac{1}{2}}}+\frac{1}{[y+2-\sqrt{y^2+4y}]^{\sigma-\frac{1}{2}}}\right\}dy,\nonumber\\
&&F'''_{\sigma}(\lambda)=e^{-2\lambda}\int_0^{\infty}\frac{(y+2)^2(3-2\lambda-\lambda y)e^{-\lambda y}}{2\sqrt{y^2+4y}}\nonumber\\
&&\quad\quad\quad\quad\quad\quad\quad\cdot\left\{\frac{1}{[y+2+\sqrt{y^2+4y}]^{\sigma-\frac{1}{2}}}+\frac{1}{[y+2-\sqrt{y^2+4y}]^{\sigma-\frac{1}{2}}}\right\}dy,
\end{eqnarray}
and
\begin{eqnarray}\label{Mar27a}
{\cal H}'_{\sigma}(x)&=&{\cal G}'_{\sigma}(e^x)e^{\frac{x}{2}}-\frac{1}{2}{\cal G}_{\sigma}(e^x)e^{-\frac{x}{2}},\nonumber\\
{\cal H}''_{\sigma}(x)&=&{\cal G}''_{\sigma}(e^x)e^{\frac{3x}{2}}+\frac{1}{4}{\cal G}_{\sigma}(e^x)e^{-\frac{x}{2}},\nonumber\\
{\cal H}'''_{\sigma}(x)&=&{\cal G}'''_{\sigma}(e^x)e^{\frac{5x}{2}}+\frac{3}{2}{\cal G}''_{\sigma}(e^x)e^{\frac{3x}{2}}+\frac{1}{4}{\cal G}'_{\sigma}(e^x)e^{\frac{x}{2}}-\frac{1}{8}{\cal G}_{\sigma}(e^x)e^{-\frac{x}{2}}.
\end{eqnarray}
Then, by (\ref{Jan28D1}) and (\ref{Janm1})--(\ref{Mar27a}), we obtain that
\begin{eqnarray}
{\cal S}_{\sigma}&=&2^{\sigma+\frac{3}{2}}\pi^{-1}\Bigg\{[\sigma^2+(1-\sigma)^2]\left[{\cal G}'_{\sigma}(1)-\frac{1}{2}{\cal G}_{\sigma}(1)\right]-\left[{\cal G}'''_{\sigma}(1)+\frac{3}{2}{\cal G}''_{\sigma}(1)+\frac{1}{4}{\cal G}'_{\sigma}(1)-\frac{1}{8}{\cal G}_{\sigma}(1)\right]\Bigg\}\nonumber\\
&=&2^{\sigma+\frac{3}{2}}\pi^{-1}\Bigg\{-{\cal G}'''_{\sigma}(1)-\frac{3}{2}{\cal G}''_{\sigma}(1)+\left[\sigma^2+(1-\sigma)^2-\frac{1}{4}\right]\left[{\cal G}'_{\sigma}(1)-\frac{1}{2}{\cal G}_{\sigma}(1)\right]\Bigg\}\nonumber\\
&=&2^{\sigma+\frac{3}{2}}\pi^{-1}\sum_{m=1}^{\infty}\sum_{n=1}^{\infty}m^{-\sigma-\frac{1}{2}}n^{\sigma-\frac{3}{2}}\Bigg\{\left[-F'''_{\sigma}(\pi mn)(\pi mn)^3-\frac{3}{2}F''_{\sigma}(\pi mn)(\pi mn)^2\right]\nonumber\\
&&\ \ \ \ +\left[\sigma^2+(1-\sigma)^2-\frac{1}{4}\right]\left[F'_{\sigma}(\pi mn)\pi mn-\frac{1}{2}F_{\sigma}(\pi mn)\right]\Bigg\}\label{Jan10b}\\
&=&2^{\sigma+\frac{1}{2}}\pi^{-1}\sum_{m=1}^{\infty}\sum_{n=1}^{\infty}m^{-\sigma-\frac{1}{2}}n^{\sigma-\frac{3}{2}}\Bigg\{\int_2^{\infty}\frac{[(\pi mny)^2-\frac{9}{2}\pi mny+3](\pi mn)^2ye^{-\pi mn y}}{\sqrt{y^2-4}}\nonumber\\
&&\ \ \ \ \ \ \ \ \ \ \ \ \ \ \ \ \ \ \ \ \ \ \ \ \ \ \ \cdot\left(\frac{1}{[y+\sqrt{y^2-4}]^{\sigma-\frac{1}{2}}}+\frac{1}{[y-\sqrt{y^2-4}]^{\sigma-\frac{1}{2}}}\right)dy\nonumber\\
&&\ \ \ \ +\left[2\left(\sigma-\frac{1}{2}\right)^2+\frac{1}{4}\right]\int_2^{\infty}\frac{(\frac{1}{2}-\pi mny)(\pi mn)e^{-\pi mn y}}{\sqrt{y^2-4}}\nonumber\\
&&\ \ \ \ \ \ \ \ \ \ \ \ \ \ \ \ \ \ \ \ \ \ \ \ \ \ \ \cdot\left(\frac{1}{[y+\sqrt{y^2-4}]^{\sigma-\frac{1}{2}}}+\frac{1}{[y-\sqrt{y^2-4}]^{\sigma-\frac{1}{2}}}\right)dy\Bigg\}\nonumber\\
&=&2\sum\limits_{m,n=1}^{\infty}\left(\frac{n}{m}\right)^{\sigma-\frac{1}{2}}\Bigg\{\int_1^{\infty}\frac{[4(\pi mnw)^2-9\pi mnw+3](2\pi mnw)e^{-2\pi mn w}}{\sqrt{w^2-1}}\nonumber\\
&&\ \ \ \ \ \ \ \ \ \ \ \ \ \ \ \ \ \ \ \ \ \ \ \ \ \ \ \cdot\left(\frac{1}{[w+\sqrt{w^2-1}]^{\sigma-\frac{1}{2}}}+\frac{1}{[w-\sqrt{w^2-1}]^{\sigma-\frac{1}{2}}}\right)dw\nonumber\\
&&\ \ \ \ +\left[2\left(\sigma-\frac{1}{2}\right)^2+\frac{1}{4}\right]\int_1^{\infty}\frac{(\frac{1}{2}-2\pi mnw)e^{-2\pi mn w}}{\sqrt{w^2-1}}\nonumber\\
&&\ \ \ \ \ \ \ \ \ \ \ \ \ \ \ \ \ \ \ \ \ \ \ \ \ \ \ \cdot\left(\frac{1}{[w+\sqrt{w^2-1}]^{\sigma-\frac{1}{2}}}+\frac{1}{[w-\sqrt{w^2-1}]^{\sigma-\frac{1}{2}}}\right)dw\Bigg\}\nonumber\\
&=&\int_1^{\infty}\Bigg\{-2y^3{\cal J}'''_{\sigma-\frac{1}{2}}(y)-9y^2{\cal J}''_{\sigma-\frac{1}{2}}(y)+y\left[4\left(\sigma-\frac{1}{2}\right)^2-\frac{11}{2}\right]{\cal J}'_{\sigma-\frac{1}{2}}(y)\nonumber\\
&&\quad\quad\quad+\left[2\left(\sigma-\frac{1}{2}\right)^2+\frac{1}{4}\right]{\cal J}_{\sigma-\frac{1}{2}}(y)\Bigg\}{\eta}_{\sigma-\frac{1}{2}}(y)dy,\nonumber
\end{eqnarray}
and
\begin{eqnarray}
{\cal T}_{\sigma}&=&2^{\sigma+\frac{3}{2}}\pi^{-1}\left[{\cal G}'_{\sigma}(1)-\frac{1}{2}{\cal G}_{\sigma}(1)\right]\nonumber\\
&=&2^{\sigma+\frac{3}{2}}\pi^{-1}\sum_{m=1}^{\infty}\sum_{n=1}^{\infty}m^{-\sigma-\frac{1}{2}}n^{\sigma-\frac{3}{2}}\left[F'_{\sigma}(\pi mn)\pi mn-\frac{1}{2}F_{\sigma}(\pi mn)\right]\label{Jan10c}\\
&=&2^{\sigma+\frac{1}{2}}\pi^{-1}\sum_{m=1}^{\infty}\sum_{n=1}^{\infty}m^{-\sigma-\frac{1}{2}}n^{\sigma-\frac{3}{2}}\int_2^{\infty}\frac{(\frac{1}{2}-\pi mny)(\pi mn)e^{-\pi mn y}}{\sqrt{y^2-4}}\nonumber\\
&&\ \ \ \ \ \ \ \ \ \ \ \ \ \ \ \ \ \ \ \ \ \ \ \ \ \ \ \cdot\left(\frac{1}{[y+\sqrt{y^2-4}]^{\sigma-\frac{1}{2}}}+\frac{1}{[y-\sqrt{y^2-4}]^{\sigma-\frac{1}{2}}}\right)dy\nonumber\\
&=&2\sum\limits_{m,n=1}^{\infty}\left(\frac{n}{m}\right)^{\sigma-\frac{1}{2}}\int_1^{\infty}\frac{(\frac{1}{2}-2\pi mnw)e^{-2\pi mn w}}{\sqrt{w^2-1}}\left(\frac{1}{[w+\sqrt{w^2-1}]^{\sigma-\frac{1}{2}}}+\frac{1}{[w-\sqrt{w^2-1}]^{\sigma-\frac{1}{2}}}\right)dw\nonumber\\
&=&\int_1^{\infty}\left[2y{\cal J}'_{\sigma-\frac{1}{2}}(y)+{\cal J}_{\sigma-\frac{1}{2}}(y)\right]{\eta}_{\sigma-\frac{1}{2}}(y)dy.\nonumber
\end{eqnarray}

By (\ref{Jan23a}), we get
\begin{eqnarray}\label{Mar28a}
F_{\sigma}(\lambda)=2^{-(\sigma+\frac{1}{2})}\int_1^{\infty}\lambda e^{-2\lambda y}\eta_{\sigma-\frac{1}{2}}(y)dy,\ \ \ \ \lambda>0.
\end{eqnarray}
Then, by (\ref{Jan28D1}),  (\ref{Jan28D2}) and (\ref{Mar28a}), we find that
\begin{eqnarray}\label{Mar26a}
&&2^{\sigma+\frac{3}{2}}\pi^{-1}\int_{0}^{\infty}{\cal H}_{\sigma}(x)\cos(tx)dx\nonumber\\
&=&2^{\sigma+\frac{3}{2}}\pi^{-1}\int_{1}^{\infty}{\cal G}_{\sigma}(\lambda)\lambda^{-\frac{3}{2}}\cos(t\ln\lambda)d\lambda\nonumber\\
&=&2\sum_{m=1}^{\infty}\sum_{n=1}^{\infty}m^{-\sigma+\frac{1}{2}}n^{\sigma-\frac{1}{2}}\int_{1}^{\infty}\int_1^{\infty} e^{-2\pi mnxy}x^{-\frac{1}{2}}\cos(t\ln x){\eta}_{\sigma-\frac{1}{2}}(y)dxdy\nonumber\\
&=&4\int_1^{\infty}\int_1^{\infty}\cos(2t\ln x){\cal J}_{\sigma-\frac{1}{2}}(x^2y){\eta}_{\sigma-\frac{1}{2}}(y)dxdy.
\end{eqnarray}
Thus, (\ref{Mar28as}) holds by (\ref{13xyzw}) and (\ref{Mar26a}).

\vskip 0.5cm
\noindent {\bf Proof of Theorem \ref{thm000} (iii).}\ \ By (i) and (ii), we obtain that
\begin{eqnarray*}
&&2|\xi(\sigma-it)|^2\nonumber\\
&=&\int_1^{\infty}\Bigg\{-2y^3{\cal J}'''_{\sigma-\frac{1}{2}}(y)-9y^2{\cal J}''_{\sigma-\frac{1}{2}}(y)+y\left[4\left(\sigma-\frac{1}{2}\right)^2-\frac{11}{2}\right]{\cal J}'_{\sigma-\frac{1}{2}}(y)\nonumber\\
&&\quad\quad\quad+\left[2\left(\sigma-\frac{1}{2}\right)^2+\frac{1}{4}\right]{\cal J}_{\sigma-\frac{1}{2}}(y)\Bigg\}{\eta}_{\sigma-\frac{1}{2}}(y)dy\nonumber\\
&&+t^2\int_1^{\infty}\left[2y{\cal J}'_{\sigma-\frac{1}{2}}(y)+{\cal J}_{\sigma-\frac{1}{2}}(y)\right]{\eta}_{\sigma-\frac{1}{2}}(y)dy\nonumber\\
&&+4\{t^2+(1-\sigma)^2\}\{t^2+\sigma^2\}\int_1^{\infty}\int_1^{\infty}\cos(2t\ln x){\cal J}_{\sigma-\frac{1}{2}}(x^2y){\eta}_{\sigma-\frac{1}{2}}(y)dxdy\\
&=&\int_1^{\infty}\Bigg\{-2y^3{\cal J}'''_{\sigma-\frac{1}{2}}(y)-9y^2{\cal J}''_{\sigma-\frac{1}{2}}(y)-\frac{[11-8(\sigma-\frac{1}{2})^2]y}{2}{\cal J}'_{\sigma-\frac{1}{2}}(y)+\frac{1+8(\sigma-\frac{1}{2})^2}{4}{\cal J}_{\sigma-\frac{1}{2}}(y)\Bigg\}{\eta}_{\sigma-\frac{1}{2}}(y)dy\\
&&+t^2\int_1^{\infty}\left[2y{\cal J}'_{\sigma-\frac{1}{2}}(y)+{\cal J}_{\sigma-\frac{1}{2}}(y)\right]{\eta}_{\sigma-\frac{1}{2}}(y)dy\\
&&+\left\{t^4+\left[2\left(\sigma-\frac{1}{2}\right)^2+\frac{1}{2}\right]t^2+\left[\left(\sigma-\frac{1}{2}\right)^2-\frac{1}{4}\right]^2\right\}\\
&&\quad\quad\cdot\sum_{k=0}^{\infty}\frac{(-1)^k2^{2(k+1)}t^{2k}}{(2k)!}\int_1^{\infty}\int_1^{\infty}(\ln x)^{2k}{\cal J}_{\sigma-\frac{1}{2}}(x^2y){\eta}_{\sigma-\frac{1}{2}}(y)dxdy.
\end{eqnarray*}
Therefore, (\ref{Feb11}) holds and the proof is complete.\hfill \fbox

\begin{rem}\label{rmk111} (i) By (\ref{RRR2}), we get
\begin{eqnarray}\label{March12a}
|\xi(\sigma-it)|^2=\frac{1}{2}\sum_{k=0}^{\infty}\Bigg\{\frac{(-1)^k}{(2k)!}\int_0^{\infty}{\cal U}_{\sigma}(y)y^{2k}dy\Bigg\}t^{2k},\ \ \ \ \forall \sigma, t\in\mathbb{R}.
\end{eqnarray}
Comparing (\ref{Feb11}) with (\ref{March12a}), we find that
$$
c_{\sigma-\frac{1}{2}}(k)=\frac{(-1)^k}{2\cdot(2k)!}\int_0^{\infty}{\cal U}_{\sigma}(y)y^{2k}dy,\ \ \ \ k\in\{0\}\cup\mathbb{N}.
$$
Then,
$$
(-1)^kc_{\sigma-\frac{1}{2}}(k)>0,\ \ \ \ k\in\{0\}\cup\mathbb{N}.
$$
Further, by (\ref{rem222}), we find that  for $\sigma\in(\frac{1}{2},1)$,
\begin{eqnarray*}
\left|c_{\sigma-\frac{1}{2}}(k)\right|&<& \frac{48\pi^8}{(2k)!}\int_0^{\infty}e^{5y-2e^y}y^{2k}dy\nonumber\\
&<&\frac{48\pi^8}{(2k)!}\left[e^{15}3^{2k+1}+\int_3^{\infty}e^{-y^2}y^{2k}dy\right]\nonumber\\
&<&\frac{48\pi^8(e^{15}3^{2k+1}+k!)}{(2k)!}.
\end{eqnarray*}

\noindent (ii) Let $\Upsilon$ be the set of all $\tau\in(0,\frac{1}{2})$ such that inequality (\ref{Jan27a}) does not hold. Then, the RH is equivalent to the statement that $\Upsilon$ equals the empty set. Since $\xi(s)$ has  countable number of distinct zeros in $\mathbb{C}$, $\Upsilon$  is at most countable.
\end{rem}

\section{Polynomial approximations}\setcounter{equation}{0}

By representation (\ref{13xyz}), we find that the RH is equivalent to the statement that for any $\sigma\in(\frac{1}{2},1)$,
\begin{eqnarray}\label{Jan29a}
{\cal S}_{\sigma}+{\cal T}_{\sigma}t^2+\{t^2+(1-\sigma)^2\}\{t^2+\sigma^2\}\int_0^{\infty}{\cal W}_{\sigma}(x)e^{-\sigma x}\cos(tx)dx>0,\ \ \ \ \forall t\in\mathbb{R}.
\end{eqnarray}
Different from some other inequalities that are equivalent to the RH, e.g., Suzuki \cite[Theorem 1.7]{S}, special arithmetic functions are not used in inequality (\ref{Jan29a}).

Let $(\Omega, {\cal F}, P)$ be a probability space. Denote by $E$ the expectation with respect to $P$. Suppose that $X_{\sigma}$ is a continuous random variable defined on $(\Omega, {\cal F}, P)$ with  probability
density function
$$\rho_{\sigma}(x):=\frac{{\cal W}_{\sigma}(x)e^{-\sigma x}1_{\{x\ge 0\}}}{\int_0^{\infty}{\cal W}_{\sigma}(y)e^{-\sigma y}dy},\ \ \ \ x\in\mathbb{R}.
$$
Then, inequality (\ref{Jan29a}) is equivalent to
\begin{eqnarray}\label{MM}
E[\cos(tX_{\sigma})]>\frac{-({\cal S}_{\sigma}+{\cal T}_{\sigma}t^2)}{\left[\int_0^{\infty}{\cal W}_{\sigma}(y)e^{-\sigma y}dy\right]\{t^2+(1-\sigma)^2\}\{t^2+\sigma^2\}},\ \ \ \ \forall t\in\mathbb{R}.
\end{eqnarray}
It is interesting to use the Monte Carlo method to numerically check inequality (\ref{MM}).

As the first step to analyze  inequality (\ref{Jan29a}), we consider polynomial approximations. For $n\in\mathbb{N}$, let $C_n$ be defined by (\ref{28A}).
\begin{lem}\label{lemb}
For $T\in\mathbb{N}$, $n\in\{0\}\cup\mathbb{N}$ and $\sigma\in(\frac{1}{2},1)$, we have
$$
\frac{T^{2n}}{(2n)!}\int_0^{\infty}{\cal W}_{\sigma}(x)e^{-\sigma x}x^{2n}dx<\frac{(C_1+C_{4{T}+3})C_{2{T}+3}T^{2n}}{[2({T}+1)]^{2n+1}}.
$$
\end{lem}

\noindent {\bf Proof.}\ \ For $x\ge 0$, we have
\begin{eqnarray*}
{\cal W}_{\sigma}(x)&\le&C_{4{T}+3}C_{2{T}+3}\int_0^{\infty}e^{-(4{T}+3)y}e^{(2{T}+3)(y-x)}e^{x+(2\sigma-1) y}dy+\int_{0}^1R(z)R(e^xz^{-1})e^xz^{2(\sigma-1)}dz\\
&<&C_{4{T}+3}C_{2{T}+3}e^{-2({T}+1)x}\int_{0}^{\infty}e^{-y}dy+C_{2{T}+3}e^{-2({T}+1)x}\int_{0}^1z^{2({T}+\sigma)+1}R(z)dz\\
&<&(C_1+C_{4{T}+3})C_{2{T}+3}e^{-2({T}+1)x}.
\end{eqnarray*}
Then,
\begin{eqnarray*}
\frac{T^{2n}}{(2n)!}\int_0^{\infty}{\cal W}_{\sigma}(x)e^{-\sigma x}x^{2n}dx&<& \frac{(C_1+C_{4{T}+3})C_{2{T}+3}T^{2n}}{(2n)!}\int_0^{\infty}e^{-2({T}+1)x}x^{2n}dx\\
&=&\frac{(C_1+C_{4{T}+3})C_{2{T}+3}T^{2n}}{[2({T}+1)]^{2n+1}}.
\end{eqnarray*}
\hfill \fbox

By the Taylor expansion, Theorem \ref{thm000}, Lemma \ref{lemb} and the fact that for any $m\in\mathbb{N}$,
$$
\sum_{n=0}^{2m-1}\frac{(-1)^n}{(2n)!}x^{2n}\le\cos x\le\sum_{n=0}^{2m-2}\frac{(-1)^n}{(2n)!}x^{2n},\ \ \ \ \forall x\in\mathbb{R},
$$
we obtain the following result.

\begin{pro}\label{thm11} (i) For any  $m\in\mathbb{N}$, any $t\in\mathbb{R}$ and any $\sigma\in(\frac{1}{2},1)$,
\begin{eqnarray*}
{\cal S}_{\sigma}+{\cal T}_{\sigma}t^2+\{t^2+(1-\sigma)^2\}\{t^2+\sigma^2\}\sum_{n=0}^{2m-2}\frac{(-1)^n}{(2n)!}t^{2n}\int_0^{\infty}{\cal W}_{\sigma}(x)e^{-\sigma x}x^{2n}dx>0.
\end{eqnarray*}

\noindent (ii) The RH is true if and only if for any $\sigma\in(\frac{1}{2},1)$ and any $T\in\mathbb{N}$, there exists $m(\sigma,T)\in\mathbb{N}$ such that for any $t\in[0,T]$,
\begin{eqnarray}\label{sos2}
{\cal S}_{\sigma}+{\cal T}_{\sigma}t^2+\{t^2+(1-\sigma)^2\}\{t^2+\sigma^2\}\sum_{n=0}^{2m(\sigma,T)-1}\frac{(-1)^n}{(2n)!}t^{2n}\int_0^{\infty}{\cal W}_{\sigma}(x)e^{-\sigma x}x^{2n}dx>0.
\end{eqnarray}
\end{pro}

\begin{rem} Note that, for any fixed $\sigma\in(\frac{1}{2},1)$,  (\ref{sos2}) is a polynomial inequality with variable $t$. It is possible to apply the polynomial sum-of-squares with the help of computer  to investigate its validity.
\end{rem}

For $\sigma\in(\frac{1}{2},1)$, $N=(N^{(1)},N^{(2)})\in \mathbb{N}^2$ and $t\ge 0$, define
$${\cal V}^N_{\sigma}(t)={\cal S}_{\sigma}+{\cal T}_{\sigma}t^2+\{t^2+(1-\sigma)^2\}\{t^2+\sigma^2\}\sum_{n=0}^{N^{(2)}}\frac{(-1)^n}{(2n)!}t^{2n}\int_0^{N^{(1)}}{\cal W}_{\sigma}(x)e^{-\sigma x}x^{2n}dx.
$$
Let ${\cal C}$ be defined by (\ref{C}). For $\varepsilon\in(0,1)$ and $T\in\mathbb{N}$, define $N_{\varepsilon,\sigma, T}=(N^{(1)}_{\varepsilon, \sigma,T},N^{(2)}_{\varepsilon, \sigma,T})$ by
$$
N^{(1)}_{\varepsilon,\sigma,T}=\left\lceil\frac{1}{\sigma}\ln\left\{\frac{8{\cal C}^2(T^2+1)^2}{\sigma\varepsilon}\right\}\right\rceil,
$$
and
$$
N^{(2)}_{\varepsilon,\sigma, T}=\left\lceil\left\{\frac{8{\cal C}^2(T^2+1)^2}{\sigma\varepsilon}\right\}^2\right\rceil.
$$
By Lemma \ref{lem22} (vii), we get
\begin{eqnarray}\label{app1}
(T^2+1)^2\int_{N^{(1)}_{\varepsilon,\sigma,T}}^{\infty}{\cal W}_{\sigma}(x)e^{-\sigma x}dx
&<&2{\cal C}^2(T^2+1)^2\int_{N^{(1)}_{\varepsilon,\sigma,T}}^{\infty}e^{-\sigma x}dx\nonumber\\
&=&\frac{2{\cal C}^2(T^2+1)^2}{\sigma}\cdot e^{-\sigma N^{(1)}_{\varepsilon,\sigma,T}}\nonumber\\
&<&\frac{\varepsilon}{4}.
\end{eqnarray}

We have the following result for the validity of the RH.
\begin{pro}\label{thmsd}
The following statements are equivalent to each other.

(i) The RH is true.

(ii) For any $\sigma\in(\frac{1}{2},1)$ and any $T\in\mathbb{N}$,
\begin{equation}\label{23c1}
\lim_{\varepsilon\rightarrow0}\frac{\min_{t\in[0,T]}{\cal V}^{N_{\varepsilon,\sigma, T}}_{\sigma}(t)}{\varepsilon}=\infty.
\end{equation}

(iii) For any $\sigma\in(\frac{1}{2},1)$ and any $T\in\mathbb{N}$,
\begin{equation}\label{23a1}
\liminf_{\varepsilon\rightarrow0}\frac{\min_{t\in[0,T]}{\cal V}^{N_{\varepsilon,\sigma, T}}_{\sigma}(t)}{\varepsilon}>0.
\end{equation}

(iv) For any $\sigma\in(\frac{1}{2},1)$ and any $T\in\mathbb{N}$, there exists $\varepsilon_{\sigma,T}\in(0,1)$ such that
\begin{equation}\label{23d3}
\min_{t\in[0,T]}{\cal V}^{N_{\varepsilon_{\sigma,T},\sigma, T}}_{\sigma}(t)\ge \frac{\varepsilon_{\sigma,T}}{2}.
\end{equation}
\end{pro}

\noindent {\bf Proof.}\ \  For $\sigma\in(\frac{1}{2},1)$ and $t\in\mathbb{R}$, define
\begin{eqnarray*}
{\cal V}_{\sigma}(t)=2|\xi(\sigma-it)|^2.
\end{eqnarray*}
Then, the RH  is true if and only if for any $t\ge0$ and
any $\sigma\in(\frac{1}{2},1)$,
$$
{\cal V}_{\sigma}(t)>0.
$$
By Lemma \ref{lem22} (vii), (\ref{13xyz}), (\ref{app1}), Taylor's theorem and Stirling's approximation, we get
\begin{eqnarray*}
&&\max_{t\in[0,T]}\left|{\cal V}_{\sigma}(t)-{\cal V}_{\sigma}^{N_{\varepsilon, \sigma,T}}(t)\right|\nonumber\\
&\le&(T^2+1)^2\left[\int_{N^{(1)}_{\varepsilon,\sigma, T}}^{\infty}{\cal W}_{\sigma}(x)e^{-\sigma x}dx+\int_0^{N^{(1)}_{\varepsilon,\sigma, T}}{\cal W}_{\sigma}(x)e^{-\sigma x}\cdot\frac{(TN^{(1)}_{\varepsilon,\sigma, T})^{2N^{(2)}_{\varepsilon,\sigma, T}+1}}{(2N^{(2)}_{\varepsilon,\sigma, T}+1)!}dx\right]\nonumber
\end{eqnarray*}
\begin{eqnarray}\label{OOO}
&<&\frac{\varepsilon}{4}+\frac{2{\cal C}^2(T^2+1)^2}{\sigma}\cdot\frac{(TN^{(1)}_{\varepsilon,\sigma, T})^{2N^{(2)}_{\varepsilon,\sigma, T}+1}}{\sqrt{2\pi(2N^{(2)}_{\varepsilon,\sigma, T}+1) }\cdot\left[e^{-1}(2N^{(2)}_{\varepsilon,\sigma, T}+1)\right]^{2N^{(2)}_{\varepsilon,\sigma, T}+1}}\nonumber\\
&<&\frac{\varepsilon}{4}+\frac{2{\cal C}^2(T^2+1)^2}{\sigma}\cdot\frac{1}{\sqrt{2\pi(2N^{(2)}_{\varepsilon,\sigma, T}+1) }}\nonumber\\
&<&\frac{\varepsilon}{2}.
\end{eqnarray}

\noindent ``(i)$\Rightarrow$(ii)":\ \ Suppose that the RH is true. Let $M$ be an arbitrary natural number. Set
\begin{eqnarray*}
\varepsilon_{\sigma,T}=\frac{\min_{t\in[0,T]}{\cal V}_{\sigma}(t)}{M}.
\end{eqnarray*}
Then, by (\ref{OOO}), we obtain that for any $\varepsilon\in(0, \varepsilon_{\sigma,T}\wedge 1)$,
\begin{eqnarray*}
\frac{\min_{t\in[0,T]}{\cal V}_{\sigma}^{N_{\varepsilon,\sigma, T}}(t)}{\varepsilon}&\ge&  \frac{\min_{t\in[0,T]}{\cal V}_{\sigma}(t)}{\varepsilon}-\frac{\max_{t\in[0,T]}\left|{\cal V}_{\sigma}(t)-{\cal V}_{\sigma}^{N_{\varepsilon, \sigma,T}}\right|}{\varepsilon}\\
&>&\frac{M\varepsilon_{\sigma,T}}{\varepsilon}-\frac{1}{2}\\
&>&M-\frac{1}{2}.
\end{eqnarray*}
Since $M$ is arbitrary, (\ref{23c1}) holds.

\noindent ``(ii)$\Rightarrow$(iii)":\ \ Obvious.

\noindent ``(iii)$\Rightarrow$(i)":\ \ Fix $\sigma\in(\frac{1}{2},1)$ and $T\in\mathbb{N}$. By (\ref{23a1}), there exist $a,\varepsilon_0\in(0,1)$ such that
\begin{eqnarray}\label{23d1}
\frac{\min_{t\in[0,T]}{\cal V}^{N_{\varepsilon,\sigma, T}}_{\sigma}(t)}{\varepsilon}>a,\ \ \ \ \forall\varepsilon\in(0,\varepsilon_0).
\end{eqnarray}
Then, by (\ref{OOO}) and (\ref{23d1}), we get
\begin{eqnarray*}
\min_{t\in[0,T]}{\cal V}_{\sigma}(t)&\ge& \min_{t\in[0,T]}{\cal V}_{\sigma}^{N_{a\varepsilon_0,\sigma, T}}(t)-\max_{t\in[0,T]}\left|{\cal V}_{\sigma}(t)-{\cal V}_{\sigma}^{N_{a\varepsilon_0, \sigma,T}}\right|\\
&>&a\varepsilon_0-\frac{a\varepsilon_0}{2}\\
&>&0.
\end{eqnarray*}
Thus, by the arbitrariness of $\sigma$ and $T$, we conclude that (iii) implies (i).

\noindent ``(ii)$\Rightarrow$(iv)":\ \ Obvious.

\noindent ``(iv)$\Rightarrow$(i)":\ \ Fix $\sigma\in(\frac{1}{2},1)$ and $T\in\mathbb{N}$. Assume that (\ref{23d3}) holds for some $\varepsilon_{\sigma,T}\in(0,1)$. Then, by (\ref{23d3}) and (\ref{OOO}), we get
\begin{eqnarray*}
\min_{t\in[0,T]}{\cal V}_{\sigma}(t)&\ge&  \min_{t\in[0,T]}{\cal V}_{\sigma}^{N_{\varepsilon_{\sigma,T},\sigma, T}}(t)-\max_{t\in[0,T]}\left|{\cal V}_{\sigma}(t)-{\cal V}_{\sigma}^{N_{\varepsilon_{\sigma,T},\sigma, T}}(t)\right|\\
&>&\frac{\varepsilon_{\sigma,T}}{2}-\frac{\varepsilon_{\sigma,T}}{2}\\
&=&0.
\end{eqnarray*}
Therefore, by the arbitrariness of $\sigma$ and $T$, we conclude that (iv) implies (i).
\hfill \fbox

\section{Equivalent description of the RH using the orthogonalization time of quantum states}\setcounter{equation}{0}

The Riemann zeta function $\zeta$ plays an integral role from
classical mechanics to statistical physics (cf. Schumayer \& Hutchinson \cite{SH}).  A lot of attempts have been made to construct suitable physical models related to $\zeta$. The Hilbert-P\'olya conjecture states that the nontrivial zeros of $\zeta$ correspond to eigenvalues of a self-adjoint operator and hence the spectral theory provides a possible approach to the RH. The Hilbert-P\'olya conjecture operator connects with quantum mechanics and has the form $\frac{1}{2}+iS$, where
$S$ is the Hamiltonian of a particle that moves under the influence of a potential
$V$. Then, the RH is equivalent to the statement that the Hamiltonian is Hermitian, or equivalently that  $V$ is real (cf. Wikipedia \cite{W2}).

Different from the existing explorations of the RH, in this section, we use a variant of inequality (\ref{Jan27a}) to connect the RH with the orthogonalization time of quantum states. Through considering a  quantum dynamics defined on the classical Wiener space and constructing suitable wave functions
by virtue of the Brownian motion, we will show that the RH is equivalent to the statement that the wave functions never evolve into orthogonal states.

\subsection{A variant inequality}

In this subsection, we present another variant of inequality (\ref{Jan27a}). By Theorem \ref{thm000}, (\ref{13xyzw}) and \cite[Proposition 2.5 (xii) and Example 2.11]{Sato}, we get
\begin{eqnarray*}
&&\frac{|\xi(\sigma-it)|^2}{\{t^2+(1-\sigma)^2\}\{t^2+\sigma^2\}}\\
&=&\frac{{\cal S}_{\sigma}+{\cal T}_{\sigma}t^2}{2\{t^2+(1-\sigma)^2\}\{t^2+\sigma^2\}}+2^{\sigma+\frac{1}{2}}\pi^{-1}\int_{0}^{\infty}{\cal H}_{\sigma}(x)\cos(tx)dx\\
&=&\frac{1}{2(2\sigma-1)}\Bigg[\frac{{\cal S}_{\sigma}-{\cal T}_{\sigma}(1-\sigma)^2}{t^2+(1-\sigma)^2}-\frac{{\cal S}_{\sigma}-{\cal T}_{\sigma}\sigma^2}{t^2+\sigma^2}\Bigg]+2^{\sigma-\frac{1}{2}}\pi^{-1}\int_{-\infty}^{\infty}{\cal H}_{\sigma}(|x|)\cos(tx)dx\\
&=&\frac{1}{4(2\sigma-1)}\int_{-\infty}^{\infty}\Bigg[\frac{{\cal S}_{\sigma}-{\cal T}_{\sigma}(1-\sigma)^2}{(1-\sigma)e^{(1-\sigma)|x|}}-\frac{{\cal S}_{\sigma}-{\cal T}_{\sigma}\sigma^2}{\sigma e^{\sigma|x|}}\Bigg]e^{-itx}dx+2^{\sigma-\frac{1}{2}}\pi^{-1}\int_{-\infty}^{\infty}{\cal H}_{\sigma}(|x|)e^{-itx}dx\\
&=&\frac{1}{4(2\sigma-1)}\int_{-\infty}^{\infty}\frac{\sigma[{\cal S}_{\sigma}-{\cal T}_{\sigma}(1-\sigma)^2]e^{(\sigma-1)|x|}-(1-\sigma)[{\cal S}_{\sigma}-{\cal T}_{\sigma}\sigma^2]e^{-\sigma|x|}}{\sigma(1-\sigma)}\cdot e^{-itx}dx\\
&&+2^{\sigma-\frac{1}{2}}\pi^{-1}\int_{-\infty}^{\infty}{\cal H}_{\sigma}(|x|)e^{-itx}dx.
\end{eqnarray*}
For $\sigma\in(\frac{1}{2},1)$, define
\begin{eqnarray}\label{Jan13c}
{\cal K}_{\sigma}(x)&=&2^{\sigma+\frac{3}{2}}\pi^{-1}\sigma(1-\sigma)(2\sigma-1){\cal H}_{\sigma}(x)\nonumber\\
&&+\sigma[{\cal S}_{\sigma}-{\cal T}_{\sigma}(1-\sigma)^2]e^{(\sigma-1)x}-(1-\sigma)[{\cal S}_{\sigma}-{\cal T}_{\sigma}\sigma^2]e^{-\sigma x},\ \ \ \ x\in\mathbb{R}.\ \
\end{eqnarray}
Then, we have the following result.

\begin{thm}\label{thm1} The RH is equivalent to the statement that for any $\sigma\in(\frac{1}{2},1)$,
\begin{eqnarray}\label{Jan18a}
\int_{-\infty}^{\infty}{\cal K}_{\sigma}(|x|)e^{-itx}dx>0,\ \ \ \ \forall t\in\mathbb{R}.
\end{eqnarray}
\end{thm}

We claim that
\begin{equation}\label{Jan22g}
{\cal K}_{\sigma}(x)>0,\ \ \ \ \forall x\ge 0.
\end{equation}
Hence, for any $\sigma\in(\frac{1}{2},1)$, ${\cal K}_{\sigma}$ is a positive, rapidly decaying and smooth function on $[0,\infty)$.

We first consider the signs of ${\cal S}_{\sigma}$ and ${\cal T}_{\sigma}$. Note that for $\lambda>3$ the signs of $F_{\sigma}(\lambda)$, $F'_{\sigma}(\lambda)$, $F''_{\sigma}(\lambda)$, $F'''_{\sigma}(\lambda)$ are $+,-,+,-$, respectively (cf. (\ref{Jan10a})). Then, by (\ref{Jan10c}), we get
$$
{\cal T}_{\sigma}<0,\ \ \ \ \sigma\in\left(\frac{1}{2},1\right).
$$
For $m,n\in\mathbb{N}$, we have
\begin{eqnarray}\label{Mar282b}
&&\left[-F'''_{\sigma}(\pi mn)(\pi mn)^3-\frac{3}{2}F''_{\sigma}(\pi mn)(\pi mn)^2\right]\nonumber\\
&&+\left[\sigma^2+(1-\sigma)^2-\frac{1}{4}\right]\left[F'_{\sigma}(\pi mn)\pi mn-\frac{1}{2}F_{\sigma}(\pi mn)\right]\nonumber\\
&>&-F'''_{\sigma}(\pi mn)(\pi mn)^3-\frac{3}{2}F''_{\sigma}(\pi mn)(\pi mn)^2+\frac{3}{4}F'_{\sigma}(\pi mn)\pi mn-\frac{3}{8}F_{\sigma}(\pi mn)\nonumber\\
&\ge&(\pi mn)^3\Bigg[-F'''_{\sigma}(\pi mn)-\frac{3}{2\pi}F''_{\sigma}(\pi mn)+\frac{3}{4\pi^2}F'_{\sigma}(\pi mn)-\frac{3}{8\pi^3}F_{\sigma}(\pi mn)\Bigg].
\end{eqnarray}
By (\ref{Jan10a}), for $\lambda>3$, we have
\begin{eqnarray}\label{Jan22f}
&&-F'''_{\sigma}(\lambda)-\frac{3}{2\pi}F''_{\sigma}(\lambda)+\frac{3}{4\pi^2}F'_{\sigma}(\lambda)-\frac{3}{8\pi^3}F_{\sigma}(\lambda)\nonumber\\
&>&e^{-2\lambda}\int_0^{\infty}\frac{(y+2)\left(\lambda-2+\frac{\lambda y}{2}\right)e^{-\lambda y}}{2\sqrt{y^2+4y}}\left\{\frac{1}{[y+2+\sqrt{y^2+4y}]^{\sigma-\frac{1}{2}}}+\frac{1}{[y+2-\sqrt{y^2+4y}]^{\sigma-\frac{1}{2}}}\right\}dy\nonumber\\
&&+\frac{3}{4\pi^2}F'_{\sigma}(\lambda)-\frac{3}{8\pi^3}F_{\sigma}(\lambda)\nonumber\\
&>&e^{-2\lambda}\int_0^{\infty}\frac{(2\lambda-4+\lambda y)e^{-\lambda y}}{2\sqrt{y^2+4y}}\left\{\frac{1}{[y+2+\sqrt{y^2+4y}]^{\sigma-\frac{1}{2}}}+\frac{1}{[y+2-\sqrt{y^2+4y}]^{\sigma-\frac{1}{2}}}\right\}dy\nonumber\\
&&+e^{-2\lambda}\int_0^{\infty}\frac{\frac{3}{8\pi^3}[2\pi-(4\pi+1)\lambda-2\pi\lambda y]e^{-\lambda y}}{2\sqrt{y^2+4y}}\nonumber\\
&&\quad\quad\quad\quad\cdot\left\{\frac{1}{[y+2+\sqrt{y^2+4y}]^{\sigma-\frac{1}{2}}}+\frac{1}{[y+2-\sqrt{y^2+4y}]^{\sigma-\frac{1}{2}}}\right\}dy\nonumber\\
&>&0.
\end{eqnarray}
Then, by (\ref{Jan10b}), (\ref{Mar282b}) and (\ref{Jan22f}), we get
$$
{\cal S}_{\sigma}>0,\ \ \ \ \sigma\in\left(\frac{1}{2},1\right).
$$

We now prove inequality (\ref{Jan22g}). By (\ref{Jan13c}), we need only show that
$$
{\cal S}_{\sigma}+\sigma(1-\sigma){\cal T}_{\sigma}>0,
$$
which is implied by
\begin{equation}\label{Jan22a}
{\cal S}_{\sigma}+\frac{1}{4}{\cal T}_{\sigma}>0.
\end{equation}
To prove (\ref{Jan22a}), by (\ref{Jan10b}) and (\ref{Jan10c}), we need only show that for $m,n\in\mathbb{N}$,
\begin{eqnarray*}
&&\left[-F'''_{\sigma}(\pi mn)(\pi mn)^3-\frac{3}{2}F''_{\sigma}(\pi mn)(\pi mn)^2\right]\nonumber\\
&&+\left[\sigma^2+(1-\sigma)^2-\frac{1}{4}\right]\left[F'_{\sigma}(\pi mn)\pi mn-\frac{1}{2}F_{\sigma}(\pi mn)\right]\\
&&+\frac{1}{4}\left[F'_{\sigma}(\pi mn)\pi mn-\frac{1}{2}F_{\sigma}(\pi mn)\right]\\
&>&0,
\end{eqnarray*}
which is implied by
\begin{eqnarray*}
-F'''_{\sigma}(\pi mn)-\frac{3}{2\pi}F''_{\sigma}(\pi mn)+\frac{1}{\pi^2}F'_{\sigma}(\pi mn)-\frac{1}{2\pi^3}F_{\sigma}(\pi mn)>0.
\end{eqnarray*}
This last inequality can be established similarly to (\ref{Jan22f}). Hence  (\ref{Jan22a}) and therefore (\ref{Jan22g}) hold.

\subsection{Wave functions and the orthogonalization time}

Motivated by  Luo et al. \cite{L}, in this subsection, we will investigate the connection between inequality (\ref{Jan18a}), equivalently the RH, and the orthogonalization time of quantum states. It is well known that the orthogonalization time of quantum states plays an important role in quantum measurement and information processing.  Time-energy uncertainty relations that bound the evolution speed, based on the  energy spread $\Delta E$ or the average energy $E$, have been discovered by Mandelstam and Tamm \cite{MT}, Margolus and Levitin \cite{ML}, and Levitin and Toffoli \cite{LT}. Their results establish the fundamental limit of the operation rates of  information processing systems.
Recently, Georgiev discussed in \cite{G}  the case that a quantum system never evolves into an
orthogonal state.

In the sequel, we use Dirac's notation of bras and kets \cite{D}. Let $\Omega=C(\mathbb{R},\mathbb{R})$, the space of all real-valued continuous functions defined on $\mathbb{R}$. We equip $\Omega$ with the locally uniform convergence topology  and denote by ${\cal B}(\Omega)$ the Borel $\sigma$-algebra. Let $P$ be the standard Wiener measure defined on $C([0,\infty),\mathbb{R})$. We extend $P$ to $(\Omega,{\cal B}(\Omega))$ by symmetry:
\begin{eqnarray*}
P(\{\omega\in \Omega:\omega(t_i)\in A_i,-n\le i\le n\})&=&
\delta_{0}(A_0)P(\{\omega\in \Omega:\omega(t_i)\in A_i,1\le i\le n\})\\
&&\cdot P(\{\omega\in \Omega:\omega(-t_i)\in A_i,-n\le i\le -1\}),
\end{eqnarray*}
where $\delta_0$ is the Dirac measure at the origin, $t_{-n}<\cdots <t_{-1}<t_0=0<t_1<\cdots < t_n$ are arbitrary real numbers, $A_i\in{\cal B}(\mathbb{R})$, $-n\le i\le n$, and $n\in\mathbb{N}$ is arbitrary. Denote by ${\cal N}$ the collection of all $P$-null sets and let ${\cal F}:=\sigma({\cal B}(\Omega)\cup{\cal N})$ be the augmentation. Set
$$
{\cal Q}=L^2(\Omega,{\cal F},P),
$$
which is a Hilbert space describing a quantum system.

For $t\in\mathbb{R}$, define
$$
{\cal F}_t=\sigma(w(s),-\infty<s\le t),
$$
and let ${\cal F}^P_t$ be the completion of ${\cal F}_t$ with respect to $P$. Denote by $E_{t}$  the orthogonal projection operator on ${\cal Q}_t:=L^2(\Omega,{\cal F}^P_t,P)$, which is a closed subspace of ${\cal Q}$. We consider the Hamiltonian $S$, which is a self-adjoint operator defined by
$$
S=\int_{-\infty}^{\infty}t dE_{t}.
$$
For $\sigma\in(\frac{1}{2},1)$, define the normalized wave function $|\psi_{\sigma}\rangle\in{\cal Q}$ by It\^o's stochastic integral with respect to the Brownian motion:
$$
|\psi_{\sigma}\rangle=\frac{\int_{-\infty}^{\infty}\sqrt{{\cal K}_{\sigma}(|x|)}d\omega(x)}{\left[2\int_0^{\infty}{\cal K}_{\sigma}(x)dx\right]^{\frac{1}{2}}}.
$$

We have
$$
e^{-it S}|\psi_{\sigma}\rangle=\frac{\int_{-\infty}^{\infty}e^{-it x}\sqrt{{\cal K}_{\sigma}(|x|)}d\omega(x)}{\left[2\int_0^{\infty}{\cal K}_{\sigma}(x)dx\right]^{\frac{1}{2}}}.
$$
Then,
$$
\langle\psi_{\sigma}|e^{-it S}|\psi_{\sigma}\rangle=\frac{\int_{0}^{\infty}{\cal K}_{\sigma}(x)\cos(tx)dx}{\int_0^{\infty}{\cal K}_{\sigma}(x)dx}.
$$
Let $\iota_{\sigma}$ be the first time that the state $|\psi_{\sigma}\rangle$	 evolves into an orthogonal state, i.e., the maximum change according to the  quantum dynamics. Mathematically, $\iota_{\sigma}$ is defined as (we put the Planck constant $\hbar= 1$)
$$
\iota_{\sigma} = \inf\left\{t > 0:\langle\psi_{\sigma}|e^{-it S}|\psi_{\sigma}\rangle=0\right\}.
$$
Therefore, by Theorem \ref{thm1}, we obtain the following result.

\begin{thm}\label{thm2} The RH is equivalent to the statement that $\iota_{\sigma} =\infty$ for any $\sigma\in(\frac{1}{2},1)$.
\end{thm}

It is interesting to point out that there are at most countably infinite number of $\sigma\in(\frac{1}{2},1)$ such that $|\psi_{\sigma}\rangle$ might evolve into an orthogonal state, since $\xi(s)$ has only countable number of distinct zeros in $\mathbb{C}$ (cf. Remark \ref{rmk111}). 
Before ending this section, we would like to pose the following question for further research.
\begin{que} Fix a $\sigma\in(\frac{1}{2},1)$, say $\sigma=\frac{3}{4}$. Can we use quantum measurement and computation to help check whether the state $|\psi_{\sigma}\rangle$ will ever evolve into a perfectly distinguishable state?
\end{que}

\section{Appendix: approximated values of ${\cal S}_{\frac{3}{4}}$ and ${\cal T}_{\frac{3}{4}}$}\setcounter{equation}{0}

For a fixed $\sigma\in(\frac{1}{2},1)$, the two constants ${\cal S}_{\sigma}$ and ${\cal T}_{\sigma}$ play an important role in inequalities (\ref{Jan29a}) and (\ref{Jan18a}) (cf. (\ref{Jan13c})). In this appendix, we use the case $\sigma=\frac{3}{4}$  as an example to show numerical calculations for ${\cal S}_{\sigma}$ and ${\cal T}_{\sigma}$.

By replacing the $\int_0^{\infty}$'s in (\ref{Jan23a}) and (\ref{Jan10a}) with $\int_{0.001}^{20}$ and replacing the $\sum_{m=1}^{\infty}\sum_{n=1}^{\infty}$'s in (\ref{Jan10b}) and (\ref{Jan10c}) with $\sum_{m=1}^{10}\sum_{n=1}^{10}$, we obtain the following approximated values using the software {\it Mathematica}:
\begin{eqnarray}\label{Jan11a}
{\cal S}_{\frac{3}{4}}=0.473929,\ \ \ \ {\cal T}_{\frac{3}{4}}=-0.0218449.
\end{eqnarray}
We may also use (\ref{Jan27a1}) and (\ref{Jan27a2}) to calculate  ${\cal S}_{\sigma}$ and ${\cal T}_{\sigma}$, respectively. By relations (\ref{for1}) and  (\ref{for2}), we get
\begin{eqnarray}\label{Jan10ab}
{\cal S}_{\sigma}&=&-[\sigma-2(1-\sigma)^2]R(1)-\sigma R'(1)+\frac{(1-\sigma)(1+2\sigma)R^2(1)}{2}-\sigma R(1)R'(1)\nonumber\\
&&+2\sigma(1-\sigma)(2\sigma-1)\int_{1}^{\infty}y^{-2\sigma}R(y)dy-2\sigma(1-\sigma)(2\sigma-1)\int_{1}^{\infty}y^{1-2\sigma}R(y)dy\nonumber\\
&&+2(1-\sigma)^2(2\sigma-1)\int_{1}^{\infty}y^{2(\sigma-1)}R(y)dy-\sigma(1-\sigma)(2\sigma-1)\int_{1}^{\infty}y^{1-2\sigma}R^2(y)dy\nonumber\\
&&+\frac{2\sigma+1}{2}\int_{1}^{\infty}y^{2\sigma}R^2(y)dy+(1-\sigma)\int_{1}^{\infty}y^{3-2\sigma}[R'(y)]^2dy\nonumber\\
&&-(1-\sigma)\int_{1}^{\infty}y^{1-2\sigma}[yR'(y)]^2dy+\sigma\int_{1}^{\infty}y^{2\sigma-1}[yR'(y)]^2dy\nonumber\\
&&-\sigma(1-\sigma)\int_{1}^{\infty}y^{2(\sigma-1)}R(y) \{(1-\sigma)[yR(y)+y-1]+[y^2R'(y)+1]\}dy\nonumber\\
&&+\sigma(1-\sigma)\int_{1}^{\infty}y^{-2\sigma}R(y^{-1})\{\sigma R(y)+yR'(y)\}dy,
\end{eqnarray}
and
\begin{eqnarray}\label{Jan10ac}
{\cal T}_{\sigma}&=&-\int_{1}^{\infty}y^{2(\sigma-1)}R(y)\{\sigma [yR(y)+y-1]+[y^2R'(y)+1]\}dy\nonumber\\
&&-\int_{1}^{\infty}z^{-2\sigma}R(z^{-1})\{\sigma R(z)-[zR'(z)+R(z)]\}dz\nonumber\\
&=&-\int_{1}^{\infty}y^{2(\sigma-1)}R(y)\{\sigma [yR(y)+y-1]+[y^2R'(y)+1]\}dy\nonumber\\
&&+\int_{1}^{\infty}y^{-2\sigma}R(y^{-1})\{(1-\sigma) R(y)+yR'(y)\}dy.
\end{eqnarray}
By replacing the $\sum_{n=1}^{\infty}$ in (\ref{Jan27b1}) with $\sum_{n=1}^{5}$ and replacing the $\int_1^{\infty}$'s in  (\ref{Jan10ab}) and (\ref{Jan10ac}) with $\int_1^{10}$, we obtain the following  approximated values using the software {\it Mathematica}:
\begin{eqnarray}\label{Jan10aee}
{\cal S}_{\frac{3}{4}}=0.38952,\ \ \ \ {\cal T}_{\frac{3}{4}}=-0.0232205.
\end{eqnarray}
We can see that the approximated values given by  (\ref{Jan11a}) and (\ref{Jan10aee}) are close. However, this second method of calculating  ${\cal S}_{\sigma}$ and ${\cal T}_{\sigma}$ is faster to implement than the previous one.  Certainly, the analysis presented here is rough and more accurate analysis could be conducted.

\end{document}